\documentclass[11pt]{amsart}
\usepackage[hmargin=32mm, vmargin=27mm, includefoot, twoside]{geometry}
\usepackage[colorlinks,bookmarksopen=true]{hyperref}
\usepackage[latin1]{inputenc}  
\usepackage[T1]{fontenc}       
\usepackage[english]{babel}
\usepackage{amssymb}
\usepackage{latexsym}
\usepackage{mathrsfs}
\usepackage{xspace}
\usepackage[all, cmtip]{xy}
\usepackage{xr}
\newcommand{\sref}[1]{\ref{structure-#1} in~\cite{Caprace-Monod_structure}}
\newtheorem{prop}{Proposition}[section]
\newtheorem{thm}[prop]{Theorem}
\newtheorem*{thm*}{Theorem}

\newtheorem*{addendum*}{Addendum}
\newtheorem{cor}[prop]{Corollary}
\newtheorem{lem}[prop]{Lemma}

\newtheorem*{convention*}{Convention}
\theoremstyle{definition}
\newtheorem*{defn*}{Definition}

\newtheorem{remark}[prop]{Remark}

\newtheorem{question}[prop]{Question}
\newtheorem*{scholium*}{Scholium}
\theoremstyle{remark}
\newtheorem{example}[prop]{Example}
\newtheorem*{example*}{Example}
\numberwithin{equation}{section}

\newcommand{\sepa}{\nobreak\vskip1.5mm\nobreak%
\begin{center}%
\vrule height 1pt width.2\hsize%
\,.\,%
\vrule height 1pt width.2\hsize%
\end{center}%
\vskip5mm}

\newcommand{\CC}{\mathbf{C}}

\newcommand{\GG}{\mathbf{G}}

\newcommand{\HH}{\mathbf{H}}

\newcommand{\LL}{\mathbf{L}}
\newcommand{\NN}{\mathbf{N}}

\newcommand{\RR}{\mathbf{R}}

\newcommand{\ZZ}{\mathbf{Z}}

\newcommand{\la}{\langle}
\newcommand{\ra}{\rangle}
\newcommand{\inv}{^{-1}}

\newcommand{\fire}{^\mathrm{(f)}}

\newcommand{\norma}{\mathscr{N}}
\newcommand{\centra}{\mathscr{Z}}
\newcommand{\Comm}{\mathrm{Comm}}
\newcommand{\QZ}{\mathscr{Q\!Z}}
\newcommand{\se}{\subseteq}
\newcommand{\ul}{\underline}
\newcommand{\wt}{\widetilde}
\newcommand{\lra}{\longrightarrow}
\newcommand{\ud}{\underline}
\newcommand{\bsl}{\backslash}
\def\bs#1.{
              \def\temp{#1}
              \ifx\temp\empty
                   \mathcal{B}
              \else
                   \mathcal{B}(#1)
              \fi
}
\newcommand{\cat}{{\upshape CAT(0)}\xspace}
\newcommand{\tangle}[2]
{\angle_\mathrm{T}(#1,#2)}
\newcommand{\aangle}[3]
{\angle_{#1}(#2,#3)}
\newcommand{\cangle}[3]
{\overline{\angle}_{#1}(#2,#3)}
\DeclareMathOperator{\chr}{char}
\DeclareMathOperator{\rank}{rank}
\DeclareMathOperator{\proj}{proj} 
  
 \DeclareMathOperator{\Ker}{Ker}

\DeclareMathOperator{\Isom}{Is}

  \DeclareMathOperator{\soc}{soc}
\newcommand{\bd}{\partial} 

\def\Aut{\mathop{\mathrm{Aut}}\nolimits}

\begin{document}
\title[Isometry groups of non-positively curved spaces: discrete subgroups]{Isometry groups of non-positively curved spaces:\\ discrete subgroups}
\author[Pierre-Emmanuel Caprace]{Pierre-Emmanuel Caprace*}
\address{UCL, 1348 Louvain-la-Neuve, Belgium}
\email{pe.caprace@uclouvain.be}
\thanks{*F.N.R.S. Research Associate}
\author[Nicolas Monod]{Nicolas Monod$^\ddagger$}
\address{EPFL, 1015 Lausanne, Switzerland}
\email{nicolas.monod@epfl.ch}
\thanks{$^\ddagger$Supported in part by the Swiss National Science Foundation}
\begin{abstract}
We study lattices in non-positively curved metric spaces. Borel density is established in that setting
as well as a form of Mostow rigidity. A converse to the flat torus theorem is provided.
Geometric arithmeticity results are obtained after a detour through superrigidity and arithmeticity of abstract lattices.
Residual finiteness of lattices is also studied. Riemannian symmetric spaces are characterised amongst \cat spaces admitting
lattices in terms of the existence of parabolic isometries.
\end{abstract}
\maketitle
\let\languagename\relax  

\section{Introduction}
Lattices in semi-simple algebraic groups have a tantalisingly rich structure; they include
arithmetic groups and more generally S-arithmetic groups over arbitrary characteristics.
The nature of these groups is shaped in part by the fact that they are realised as isometries of a canonical
non-positively curved space: the associated Riemannian symmetric space, or Bruhat--Tits building, or a product
of both types.

\smallskip

Many other groups of rather diverse origins share this property to occur as lattices in non-positively curved
spaces, singular or not:

\begin{itemize}
\item[---] The fundamental group of a closed Riemannian manifold of non-positive sectionnal curvature. Here
the space acted upon is the universal covering, which is a Hadamard manifold.

\item[---] Many Gromov-hyperbolic groups admit a properly discontinuous cocompact action on some
{\upshape CAT($-1$)}\xspace space by isometries.
Amongst the examples arising in this way are hyperbolic Coxeter groups~\cite{Moussong},
$C'(\frac{1}{6})$ and $C'(\frac{1}{4})$-$T(4)$ small cancellation groups~\cite{Wise}, $2$-dimensional
$7$-systolic groups~\cite{JanuszSwiat}. It is in fact a well known open problem of M.~Gromov
to construct an example of a Gromov-hyperbolic group which is \emph{not} a \cat group
(see~\cite[7.B]{Gro93};  also Remark~2.3(2) in Chapter III.$\Gamma$ of~\cite{Bridson-Haefliger}).

\item[---] In~\cite{Burger-Mozes2}, striking examples of finitely presented simple groups are
constructed as lattices in a product of two locally finite trees. Tree lattices were previously
studied in~\cite{BassLubotzky}. 

\item[---] A minimal adjoint Kac--Moody group over a finite field, as defined by J.~Tits~\cite{Tits87}, is
endowed with two $BN$-pairs which yield strongly transitive actions on a pair of twinned buildings.
When the order of the ground field is large enough, the Kac--Moody group is a lattice in the product of these
two buildings~\cite{RemCRAS}.
\end{itemize}

\medskip

Subsuming all the above examples, we define a \textbf{\cat lattice}\index{CAT(0) lattice@\cat lattice}
as a pair $(\Gamma, X)$ consisting of a proper \cat space $X$ with cocompact isometry group $\Isom(X)$ and a lattice
subgroup $\Gamma < \Isom(X)$, \emph{i.e.} a discrete subgroup of finite invariant covolume (the compact-open
topology makes $\Isom(X)$ a locally compact second countable group which is thus canonically endowed with Haar
measures). We say that $(\Gamma, X)$ is \textbf{uniform}\index{CAT(0) lattice@\cat lattice!uniform} if $\Gamma$
is cocompact in $\Isom(X)$ or, equivalently, if the quotient $\Gamma \backslash X$ is compact; that case
corresponds to $\Gamma$ being a \textbf{\cat group}\index{CAT(0) group@\cat group} in the usual terminology.

\medskip

Amongst \cat lattices, the most important, and also the best understood, notably through the work of
G.~Margulis, consist undoubtedly of those arising from lattices in semi-simple groups over local fields. It is
therefore natural to address two sets of questions.

\begin{itemize}
\item[(a)] \emph{What properties of these lattices are shared by all \cat lattices?}

\item[(b)] \emph{What properties characterise them within the class of \cat lattices?}
\end{itemize}

\smallskip

This article is devoted to the study of \cat lattices and centres largely around the above questions, though
we also address the general question of the interplay between the algebraic structure of a \cat lattice and
the geometric properties of the underlying space. Some of the techniques established in the present paper have been used 
in a subsequent investigation of lattices in products of Kac--Moody groups~\cite{Caprace-Monod_KM}.

\bigskip

We shall now describe the main results of this article; for many of them, the core of the text will contain a stronger,
more precise but perhaps more ponderous version.
Our notation is standard, as recalled in the Notation section of the companion paper~\cite{Caprace-Monod_structure}.
We refer to the latter for terminology and shall quote it freely.

\sepa

\noindent\textbf{Geometric Borel density.}
As a link between the general theory exposed in~\cite{Caprace-Monod_structure} and the study of \cat
lattices, we propose the following analogue of A.~Borel's
density theorem~\cite{Borel60}\index{Borel density}.

\begin{thm}\label{thm:density-intro}
Let $X$ be a proper \cat space, $G$ a locally compact group acting continuously by isometries on $X$ and
$\Gamma < G$ a lattice. Suppose that $X$ has no Euclidean factor. 

If $G$ acts minimally without fixed point at infinity, so does $\Gamma$.
\end{thm}

\noindent This conclusion fails for spaces with a Euclidean factor. The theorem will be established more
generally for closed subgroups with finite invariant covolume.  It should be compared to (and can of course be gainfully combined
with) a similar density property of normal subgroups established as Theorem~\sref{thm:geometric_simplicity}.

\begin{remark}
Theorem~\ref{thm:density-intro} applies to general proper \cat spaces.
It implies in particular the classical Borel density theorem (see the end of Section~\ref{sec:density}).
\end{remark}

As with classical Borel density, we shall use this density statement to derive statements about the centraliser, normaliser
and radical of lattices in Section~\ref{sec:density}.

\smallskip

A more elementary variant of the above theorem shows that a large class of groups have rather restricted actions
on proper \cat spaces; as an application, one shows:

\smallskip
\itshape
Any isometric action of R.~Thompson's group $F$ on any proper \cat
space $X$ has a fixed point in $\overline{X}$\upshape,  

\smallskip
\noindent
see Corollary~\ref{cor:Thompson}. Theorem~\ref{thm:density-intro} also provides additional information about the
totally disconnected groups $D_j$ occurring in Theorem~\sref{thm:Decomposition}.
\sepa

\noindent\textbf{Lattices: Euclidean factor, boundary, irreducibility and Mostow rigidity.}
Recall that the \emph{Flat Torus
theorem}, originating in the work of Gromoll--Wolf~\cite{Gromoll-Wolf} and Lawson--Yau~\cite{Lawson-Yau},
associates Euclidean subspaces~$\RR^n$ to any subgroup~$\ZZ^n$ of a \cat group, see~\cite[\S\,II.7]{Bridson-Haefliger}.
(In the classical setting, when the \cat group is given by a compact non-positively
curved manifold, this amounts to the seemingly more symmetric statement that such a subgroup exists if and only
if there is a flat torus is the manifold.)

The converse is a well known open problem stated by M.~Gromov in~\cite[\S\,$6.\mathrm{B}_3$]{Gro93}; for manifolds see S.-T.~Yau,
problem~65 in~\cite{YauPB}). Point~\eqref{pt:fglattices:flatTT} in the following result is a (very partial) answer; in
the special case of cocompact Riemannian manifolds, this was the main result of P.~Eberlein's article~\cite{Eberlein83}.

\begin{thm}\label{thm:fglattices}
Let $X$ be a proper \cat space, $G<\Isom(X)$ a closed subgroup acting minimally and cocompactly
on $X$ and $\Gamma < G$ a finitely generated lattice. Then:
\begin{enumerate}
\item If the Euclidean factor of $X$ has dimension $n$, then $\Gamma$ possesses a finite index subgroup $\Gamma_0$
which splits as $\Gamma_0 \simeq \ZZ^n \times \Gamma'$.\label{pt:fglattices:flatTT}

\noindent Moreover, the dimension of the Euclidean factor is characterised as the maximal rank of a free Abelian
normal subgroup of $\Gamma$.

\item $G$ has no fixed point at infinity; the set of $\Gamma$-fixed points at infinity is contained in the
(possibly empty) boundary of the Euclidean factor.\label{pt:fglattices:infty}
\end{enumerate}
\end{thm}

Point~\eqref{pt:fglattices:infty} is particularly useful in conjuction with the many results assuming the absence of
fixed points at infinity in~\cite{Caprace-Monod_structure}.
In addition, it is already a first indication that the mere existence of a (finitely
generated) lattice is a serious restriction on a proper \cat space even within the class of cocompact minimal
spaces. We recall that E.~Heintze~\cite{Heintze74} produced simply connected negatively curved Riemannian
manifolds that are homogeneous (in particular, cocompact) but have a point at infinity fixed by all isometries.

\medskip

Since a \cat lattice consists of a group and a space, there are two natural notions of irreducibility: of
the group or of the space. In the case of lattices in semi-simple groups, the two notions are known to coincide
by a result of Margulis~\cite[II.6.7]{Margulis}. We prove that this is the case for \cat lattices as above.

\begin{thm}\label{thm:irred-intro}
In the setting of Theorem~\ref{thm:fglattices},
$\Gamma$ is irreducible as an abstract group if and only if for any finite index subgroup $\Gamma_1$ and
any $\Gamma_1$-equivariant decomposition $X = X_1 \times X_2$ with $X_i$ non-compact, the projection of
$\Gamma_1$ to both $\Isom(X_i)$ is non-discrete.
\end{thm}

The combination of Theorem~\ref{thm:irred-intro}, Theorem~\ref{thm:fglattices} and of an appropriate form
of superrigidity allow us to give a \cat version of Mostow rigidity for reducible spaces (Section~\ref{sec:Mostow}).

\sepa

\noindent\textbf{Geometric arithmeticity.}
We now expose results giving perhaps unexpectedly strong conclusions for \cat lattices~--- both for the group and
for the space. These results were announced in~\cite{Caprace-MonodCRAS}
in the case of \cat groups; the present setting of finitely generated lattices is more general since \cat groups are
finitely generated (\emph{cf.} Lemma~\ref{lem:finiteGeneration} below).

We recall that an isometry $g$ is \textbf{parabolic}\index{parabolic|see{isometry}}\index{isometry!parabolic} if
the translation length\index{translation length} $\inf_{x\in X}d(gx, x)$ is not achieved. For general \cat
spaces, parabolic isometries are not well understood; in fact, ruling out their existence can sometimes be the
essential difficulty in rigidity statements.

\begin{thm}\label{thm:arith:para-intro}
Let $(\Gamma, X)$ be an irreducible finitely generated \cat lattice with $X$ geodesically complete. Assume that
$X$ possesses some parabolic isometry.

If $\Gamma$ is residually finite, then $X$ is a product of symmetric spaces and Bruhat--Tits buildings. In
particular, $\Gamma$ is an arithmetic lattice unless $X$ is a real or complex hyperbolic space.

If $\Gamma$ is not residually finite, then $X$ still splits off a symmetric space factor. Moreover, the finite
residual $\Gamma_D$ of $\Gamma$ is infinitely generated and $\Gamma/\Gamma_D$ is an arithmetic group.
\end{thm}

(Recall that the \textbf{finite residual}\index{finite residual|see{residual}}\index{residual!finite} of a group
is the intersection of all finite index subgroups.)

We single out a purely geometric consequence.

\begin{cor}\label{cor:para-intro}
Let $(\Gamma, X)$ be a finitely generated \cat lattice with $X$ geodesically complete.

Then $X$ possesses a parabolic isometry if and only if  $X \cong M \times X'$, where $M$ is a symmetric space of non-compact type.
\end{cor}

Without the assumption of geodesic completeness, we still obtain an arithmeticity statement when the underlying space admits some parabolic
isometry that is \textbf{neutral}\index{neutral|see{isometry}}\index{isometry!neutral parabolic}, \emph{i.e.}
whose displacement length vanishes. Neutral parabolic isometries are even less understood, not even for their
dynamical properties (which can be completely wild at least in Hilbert space~\cite{Edelstein}); as for familiar
examples, they are provided by unipotent elements in semi-simple algebraic groups.

\begin{thm}\label{thm:Arithmeticity:neutral-intro}
Let $(\Gamma, X)$ be an irreducible finitely generated \cat lattice. If $X$ admits any neutral parabolic isometry, then either:
\begin{enumerate}
\item $\Isom(X)$ is a rank one simple Lie group with trivial centre; or:

\item $\Gamma$ has a normal subgroup $\Gamma_D$ such that $\Gamma/\Gamma_D$ is an arithmetic group. Moreover,
$\Gamma_D$ is either finite or infinitely generated.
\end{enumerate}
\end{thm}

We turn to another type of statement of arithmeticity/geometric superrigidity. Having established an abstract arithmeticity theorem
(presented below as Theorem~\ref{thm:arith:general:intro}), we can appeal to our geometric results and prove the following.

\begin{thm}\label{thm:arith:geometric:lin-intro}
Let $(\Gamma, X)$ be an irreducible finitely generated \cat lattice with $X$ geodesically complete. Assume that $\Gamma$ possesses
some faithful finite-dimensional linear representation $($in characteristic~$\neq 2,3$$)$.

If $X$ is reducible, then $\Gamma$ is an arithmetic lattice and $X$ is a product of symmetric spaces and
Bruhat--Tits buildings.
\end{thm}

Section~\ref{sec:geometric:arith} contains more results of this nature but also demonstrates by a family of examples
that some of the intricacies in the more detailed statements reflect indeed the existence of more exotic pairs
$(\Gamma, X)$.

\sepa

\noindent\textbf{Abstract arithmeticity.}
When preparing for the proof of our geometric arithmeticity statements, we are led to 
study irreducible lattices in products of general topological groups in the abstract.
Building notably on ideas of Margulis, we stablish the following arithmeticity statement
(for which we recall that the \textbf{quasi-centre}\index{quasi-centre} $\QZ$ of
a topological group is the subset of elements with open centraliser).

\begin{thm}\label{thm:arith:general:intro}
Let $\Gamma<G=G_1\times\cdots\times G_n$ be an irreducible finitely generated lattice, where each $G_i$ is
any locally compact group.

If $\Gamma$ admits a faithful Zariski-dense representation in a semi-simple group over some field of
characteristic~$\neq 2,3$, then the amenable radical $R$ of $G$ is compact and the quasi-centre $\QZ(G)$ is
virtually contained in $\Gamma\cdot R$. Furthermore, upon replacing $G$ by a finite index subgroup, the quotient
$G/R$ splits as $G^+ \times \QZ(G/R)$ where $G^+$ is a semi-simple algebraic group and the image of $\Gamma$ in
$G^+$ is an arithmetic lattice.
\end{thm}

In particular, the quasi-centre $\QZ(G/R)$ is \emph{discrete}.
In shorter terms, this theorem states that up to a compact extension, $G$ is the direct product of a semi-simple algebraic
group by a (possibly trivial) discrete group, and that the image of $\Gamma$ in the non-discrete part is an arithmetic group.
The assumption on the characteristic can be slightly weakened.

\smallskip

In the course of the proof, we characterise all irreducible finitely generated lattices in products of the form
$G=S\times D$ where $S$ is a semi-simple Lie group and $D$ a totally disconnected group
(Theorem~\ref{thm:arith:SxD}). In particular, it turns our that $D$ must necessarily be locally profinite by
analytic. The corresponding question for simple algebraic groups instead of Lie groups is also investigated
(Theorem~\ref{thm:MixedProducts}).

\sepa

\noindent\textbf{Unique geodesic extension.}
Complete simply connected Riemannian manifolds of non-positive curvature, sometimes also called Hadamard
manifolds, form a classical family of proper \cat spaces to which the preceding results may be applied. In fact,
the natural class to consider in our context consists of those proper \cat spaces in which every geodesic
segment extends \emph{uniquely} to  a bi-infinite geodesic line. Clearly,  this class contains all Hadamard
manifolds, but it presumably contains more examples. It is, however, somewhat restricted with respect to the main thrust
of the present work since it does not allow for, say, simplicial complexes; accordingly, the conclusions
of the theorem below are also more stringent.

\begin{thm}\label{thm:Hadamard}
Let $X$ be a proper \cat space with uniquely extensible geodesics. Assume that $\Isom(X)$ acts cocompactly
without fixed points at infinity.
\begin{enumerate}
\item If $X$ is irreducible, then either $X$ is a symmetric space or $\Isom(X)$ is discrete.

\item If $\Isom(X)$ possesses a finitely generated non-uniform lattice $\Gamma$ which is irreducible as an
abstract group, then $X$ is a symmetric space (without Euclidean factor).

\item Suppose that $\Isom(X)$ possesses a finitely generated lattice $\Gamma$ (if $\Gamma$ is uniform, this is
equivalent to the condition that $\Gamma$ is a discrete cocompact group of isometries of $X$). If $\Gamma$ is
irreducible (as an abstract group) and $X$ is reducible, then $X$ is a symmetric space (without Euclidean
factor).
\end{enumerate}
\end{thm}

In the special case of Hadamard manifolds, statement~(i) was known under the assumption that $\Isom(X)$
satisfies the duality condition (without assuming that $\Isom(X)$ acts cocompactly without fixed points at
infinity). This is due to P.~Eberlein (Proposition~4.8 in~\cite{Eberlein82}). Likewise,
statement~(iii) for manifolds is Proposition~4.5 in~\cite{Eberlein82}.

More recently, Farb--Weinberger~\cite{Farb-Weinberger} investigated analogous questions for aspherical manifolds.

\sepa

\noindent\textbf{Lattices and the de Rham decomposition.}
In~\cite{Caprace-Monod_structure}, we proved a ``de Rham'' decomposition
\begin{equation}\label{eq:deRham:intro}
X' \cong\ X_1\times \cdots \times X_p \times \RR^n \times Y_1\times \cdots \times Y_q
\end{equation}
for proper \cat spaces $X$ with finite-dimensional Tits boundary and such that $\Isom(X)$ has no fixed point at infinity,
see Addendum~\sref{addendum}. (Here $X'\se X$ is the canonical minimal invariant subspace, and we recall that $X'=X$
\emph{e.g.} when $X$ is geodesically complete and admits a cocompact lattice by Lemma~\sref{lem:cocompact:minimal}.)

\smallskip

It turns out that this de Rham decomposition is an invariant of \cat groups in the following sense
(see Corollary~\ref{cor:deRham:lattice}).

\begin{thm}
Let $X$ be a proper \cat space and $\Gamma < \Isom(X)$ be a group acting properly discontinuously and
cocompactly.

Then any other such space admitting a proper cocompact $\Gamma$-action has the same number of factors
in~\eqref{eq:deRham:intro} and the Euclidean factor has same dimension.
\end{thm}

\newpage\tableofcontents\newpage

\section{An analogue of Borel density}\label{sec:density}
Before discussing our analogue of Borel's density theorem~\cite{Borel60} in Section~\ref{sec:borel} below,
we present a more elementary phenomenon based on co-amenability.

\subsection{Fixed points at infinity}
Recall that a subgroup $H$ of a topological group $G$ is \textbf{co-amenable}\index{co-amenable} if any
continuous affine $G$-action on a convex compact set (in a Hausdorff locally convex topological vector space)
has a fixed point whenever it has an $H$-fixed point. The arguments of Adams--Ballmann~\cite{AB98} imply the
following preliminary step towards Theorem~\ref{thm:density}:

\begin{prop}\label{prop:coAmenable}
Let $G$ be a topological group with a continuous isometric action on a proper \cat space $X$ without Euclidean
factor. Assume that the $G$-action is minimal and does not have a global fixed point in $\bd X$.

Then any co-amenable subgroup of $G$ still has no global fixed point in $\bd X$.
\end{prop}

\begin{proof}
Suppose for a contradiction that a co-amenable subgroup $H<G$ fixes $\xi\in\bd X$. Then $G$ preserves a
probability measure $\mu$ on $\bd X$ and we obtain a convex function $f: X\to\RR$ by integrating Busemann
functions against this measure; as in~\cite{AB98}, the cocycle equation for Busemann functions (see
\S\,\sref{sec:notation})
imply that $f$ is $G$-invariant up to constants. The arguments therein show that $f$ is constant and that $\mu$
is supported on flat points. However, in the absence of a Euclidean factor, the set of flat points has a unique
circumcentre when non-empty~\cite[1.7]{AB98}; this provides a $G$-fixed point, a contradiction.
\end{proof}

Combining the above with the splitting methods used in Theorem~\sref{thm:dichotomy}, we record a consequence
showing that the exact conclusions of the Adams--Ballmann theorem~\cite{AB98}
hold under much weaker assumptions than the amenability of $G$.

\begin{cor}\label{cor:commuting:coamenable}
Let $G$ be a topological group with a continuous isometric action on a proper \cat space $X$.
Assume that $G$ contains two commuting co-amenable subgroups.

Then either $G$ fixes a point at infinity or it preserves a Euclidean subspace in $X$.
\end{cor}

We emphasise that one can easily construct a wealth of examples of highly non-amenable groups satisfying these
assumptions. For instance, given \emph{any} group $Q$, the restricted wreath product $G = \ZZ \ltimes \bigoplus_{n\in \ZZ} Q$
contains the pair of commuting co-amenable groups $H_+= \bigoplus_{n\geq 0} Q$ and $H_-= \bigoplus_{n<0} Q$, see~\cite{Monod-Popa}.
(In fact, one can even arrange for $H_\pm$ to be conjugated upon replacing $\ZZ$ by the infinite dihedral group.)

\smallskip

For similar reasons, we deduce the following fixed-point property for R.~Thompson's group\index{Thompson's group}
$$F\ :=\ \big\langle\, g_i, i\in \NN\ |\ g_i\inv g_j g_i = g_{j+1}\ \forall\,j>i \,\big\rangle;$$
this fixed-point result explains why the strategy proposed in~\cite{Farley08} to disprove amenability of $F$
with the Adams--Ballmann theorem cannot work.

\begin{cor}\label{cor:Thompson}
Any $F$-action by isometries on any proper \cat space $X$ has a fixed point in $\overline{X}$.
\end{cor}

\begin{proof}[Proof of Corollary~\ref{cor:commuting:coamenable}]
We assume that $G$ has no fixed point at infinity. By Proposition~\sref{prop:EasyDichotomy}, there is a minimal non-empty closed convex
$G$-invariant subspace. Upon considering the Euclidean decomposition~\cite[II.6.15]{Bridson-Haefliger} of the latter, we can assume
that $X$ is $G$-minimal and without Euclidean factor and need to show that $G$ fixes a point in $X$.

Let $H_\pm< G$ be the commuting co-amenable groups. In view of Proposition~\ref{prop:coAmenable}, both act without fixed point at infinity.
In particular, we have an action of $H=H_+\times H_-$ without fixed point at infinity and the splitting theorem from~\cite{Monod_superrigid}
provides us with a canonical subspace $X_+\times X_-\se X$ with component-wise and minimal $H$-action. All of $\bd X_+$ is fixed by $H_-$, which
means that this boundary is empty. Since $X$ is proper, it follows that $X_+$ is bounded and hence reduced to a point by minimality.
Thus $H_+$ fixes a point in $X \se \overline{X}$ and co-amenability implies that $G$ fixes a probability measure $\mu$ on $\overline{X}$.
If $\mu$ were supported on $\bd X$, the proof of Proposition~\ref{prop:coAmenable} would provide a $G$-fixed point at infinity, which is absurd.
Therefore $\mu(X)>0$. Now choose a bounded set $B\se X$ large enough so that $\mu(B)>\mu(X)/2$. Then any $G$-translate of $B$ must meet $B$. It follows
that $G$ has a bounded orbit and hence a fixed point as claimed.
\end{proof}

\begin{proof}[Proof of Corollary~\ref{cor:Thompson}]
We refer to~\cite{Cannon-Floyd-Parry} for a detailed introduction to the group $F$. In particular, $F$ can be realised as the group of all
orientation-preserving piecewise affine homeomorphisms of the interval $[0,1]$ that have dyadic breakpoints and slopes~$2^n$ with $n\in \ZZ$.
Given a subset $A\se[0, 1]$ we denote by $F_A<F$ the subgroup supported on $A$.

We claim that whenever $A$ has non-empty interior, $F_A$ is co-amenable in $F$. The argument is analogous to~\cite{Monod-Popa} and
to~\cite[\S\,4.F]{Glasner-Monod}; indeed, in view of the alternative definition of $F$ just recalled,
one can choose a sequence $\{g_n\}$ in $F$ such that $g_n A$ contains $[1/n, 1-1/n]$ and thus
$F_A^{g_n}$ contains $F_{[1/n, 1-1/n]}$. Consider the compact space of means
on $F/F_A$, namely finitely additive measures, endowed with the weak-* topology from the dual of $\ell^\infty(F/F_A)$.
Any accumulation point $\mu$ of the sequence of Dirac masses at $g_n\inv F_A$ will be invariant under the union $F'$ of the
groups $F_{[1/n, 1-1/n]}$. Now $F'$ is the kernel of the derivative homomorphism $F\to 2^\ZZ\times 2^\ZZ$ at the pair
of points $\{0,1\}$. In particular, $F'$ is co-amenable in $F$ and thus the $F'$-invariance of $\mu$ implies that there is also a $F$-invariant
mean on $F/F_A$, which is one of the characterisations of co-amenability~\cite{Eymard72}.

Let now $X$ be any proper \cat space with an $F$-action by isometries. We can assume that $F$ has no fixed point at infinity and therefore we can
also assume that $X$ is minimal by Proposition~\sref{prop:EasyDichotomy}. The above claim provides us with many pairs of commuting co-amenable
subgroups upon taking disjoint sets of non-empty interior. Therefore, Corollary~\ref{cor:commuting:coamenable} shows that $X\cong\RR^n$ for some $n$.
In particular the isometry group is linear. Since $F$ is finitely generated (by $g_0$ and $g_1$ in the above presentation, compare
also~\cite{Cannon-Floyd-Parry}), Malcev's theorem~\cite{Malcev40} implies that the image of $F$ is residually finite. The derived subgroup of $F$
(which incidentally coincides with the group $F'$ introduced above) being simple~\cite{Cannon-Floyd-Parry}, it follows that it acts trivially.
It remains only to observe that two commuting isometries of $\RR^n$ always have a common fixed point in $\overline{\RR^n}$, which
is a matter of linear algebra.
\end{proof}

The above reasoning can be adapted to yield similar results for branch groups and related groups; we shall address these questions elsewhere.

\subsection{Geometric density for subgroups of finite covolume}\label{sec:borel}
The following geometric density theorem generalises Borel's density (see Proposition~\ref{prop:Karpelevich-Mostow} below)
and contains Theorem~\ref{thm:density-intro} from the Introduction.

\begin{thm}\label{thm:density}\index{Borel density}
Let $G$ be a locally compact group with a continuous isometric  action on a proper \cat space $X$ without
Euclidean factor.

If $G$ acts minimally and without global fixed point in $\bd X$, then any closed subgroup with finite
invariant covolume in $G$ still has these properties.
\end{thm}

\begin{remark}
For a related statement without the assumption on the Euclidean factor of $X$ or on fixed points at infinity, see
Theorem~\ref{thm:Lattice=>NoFixedPointAtInfty_bis} below.
\end{remark}

\begin{proof}
Retain the notation of the theorem and let $\Gamma<G$ be a closed subgroup of finite invariant covolume. In
particular, $\Gamma$ is co-amenable and thus has no fixed points at infinity by
Proposition~\ref{prop:coAmenable}. By Proposition~\sref{prop:EasyDichotomy}, there is a minimal non-empty closed
convex $\Gamma$-invariant subset $Y\se X$ and it remains to show $Y=X$. Choose a point $x_0\in X$ and define
$f:X\to \RR$ by
$$f(x) = \int_{G/\Gamma} \big( d(x, gY) - d(x_0, gY)\big)\,dg.$$
This integral converges because the integrand is bounded by $d(x, x_0)$. The function $f$ is continuous, convex
(by~\cite[II.2.5(1)]{Bridson-Haefliger}) and  ``quasi-invariant'' in the sense that it satisfies
\begin{equation}\label{eq:cocycle}
f(h x) = f(x) - f(h x_0)\kern10mm \forall\,h\in G.
\end{equation}
Since $G$ acts minimally and without fixed point at infinity, this implies that $f$ is constant (see Section~2
in~\cite{AB98}; alternatively, when $\bd X$ is finite-dimensional, it follows from
Theorem~\sref{thm:geometric_simplicity} since~\eqref{eq:cocycle} implies that $f$ is invariant under the derived
subgroup $G'$).

In particular, $d(x, gY)$ is affine for all $g$. It follows that for all $x\in X$ the closed set
$$Y_x\ =\ \big\{z\in X : d(z, Y) = d(x, Y)\big\}$$
is convex. We claim that it is parallel to $Y$ in the sense that $d(z, Y)= d(y, Y_x)$ for all $z\in Y_x$ and all
$y\in Y$. Indeed, on the one hand $d(z, Y)$ is constant over $z\in Y_x$ by definition, and on the other hand
$d(y, Y_x)$ is constant by minimality of $Y$ since $d(\cdot, Y_x)$ is a convex $\Gamma$-invariant function. In
particular, $Y_x$ is $\Gamma$-equivariantly isometric to $Y$ \emph{via} nearest point projection
(compare~\cite[II.2.12]{Bridson-Haefliger}) and each $Y_x$ is $\Gamma$-minimal. At this point, Remarks~39
in~\cite{Monod_superrigid} show that there is an isometric $\Gamma$-invariant splitting
$$X\ \cong\ Y\times T.$$
It remains to show that the ``space of components'' $T$ is reduced to a point. Let thus $s,t\in T$ and let $m$
be their midpoint. Applying the above reasoning to the choice of minimal set $Y_0$ corresponding to $Y\times
\{m\}$, we deduce again that the distance to $Y_0$ is an affine function on $X$. However, this function is
precisely the distance function $d(\cdot, m)$ in $T$ composed with the projection $X\to T$. Being non-negative
and affine on $[s, t]$, it vanishes on that segment and hence $s=t$.
\end{proof}

\begin{remark}
When $\Gamma$ is cocompact in $G$, the proof can be shortened by integrating just $d(x, gY)$ in the definition
of $f$ above.
\end{remark}

\begin{cor}\label{cor:LatticeNormaliser}
Let $X$ be a proper \cat space without Euclidean factor such that $G = \Isom(X)$ acts  minimally without fixed
point at infinity, and let $\Gamma < G$ be a closed subgroup with finite invariant covolume. Then:
\begin{enumerate}
\item $\Gamma$ has trivial amenable radical.\label{pt:LatticeNormaliser:ramen}

\item The centraliser $\centra_G(\Gamma)$ is trivial.\label{pt:LatticeNormaliser:centra}

\item If $\Gamma$ is finitely generated, then is has finite index in its normaliser $\norma_G(\Gamma)$
and the latter is a finitely generated lattice in $G$.
\end{enumerate}
\end{cor}

\begin{proof}
(i) and (ii) follow by the same argument as in the proof of Theorem~\sref{thm:geometric_simplicity}. For~(iii) we
follow~\cite[Lemma~II.6.3]{Margulis}. Since $\Gamma$ is closed and countable, it is discrete by Baire's category theorem
and thus is a lattice in $G$. Since it is finitely generated, its automorphism group is countable.
By~(ii), the normaliser $\norma_G(\Gamma)$ maps injectively to $\Aut(\Gamma)$ and hence is countable as well.
Thus $\norma_G(\Gamma)$, being closed in $G$, is discrete by applying Baire again.
Since it contains the lattice $\Gamma$, it is itself a lattice and the index of
$\Gamma$ in $\norma_G(\Gamma)$ is finite. Thus $\norma_G(\Gamma)$ is finitely generated.
\end{proof}

As pointed out by P.~de la~Harpe, point\eqref{pt:LatticeNormaliser:centra} implies in particular that
any lattice in $G$ is ICC (which means by definition that all its non-trivial conjugacy classes are infinite).
As is well known, this is the criterion ensuring that the type $\mathrm{II}_1$
von Neumann algebra\index{factor!von Neumann algebra}
associated to the lattice is a factor~\cite[\S\,V.7]{Takesaki}.

\medskip

Finally, we indicate why Theorem~\ref{thm:density} implies the classical Borel density theorem
of~\cite{Borel60}. It suffices to justify the following:

\begin{prop}\label{prop:Karpelevich-Mostow}
Let $k$ by a local field (Archimedean or not), $\GG$ a semi-simple $k$-group without $k$-anisotropic factors,
$X$ the symmetric space or Bruhat--Tits building associated to $G=\GG(k)$ and $L<G$ any subgroup. If the
$L$-action on $X$ is minimal without fixed point at infinity, then $L$ is Zariski-dense.
\end{prop}

\begin{proof}
Let $\bar L$ be the ($k$-points of the) Zariski closure of $L$. Then $\bar L$ is semi-simple; this follows
\emph{e.g.} from a very special case of Corollary~\sref{cor:NoEuclideanFactor}, which guarantees that the radical
of $\bar L$ is trivial.

In the Archimedean case, we may appeal to Karpelevich--Mostow theorem (see~\cite{Karpelevich} or~\cite{Mostow}):
any semi-simple subgroup has a totally geodesic orbit in the symmetric space. So the only semi-simple subgroup
acting minimally is $G$ itself.

In the non-Archimedean case, we could appeal to E.~Landvogt functoriality theorem~\cite{Landvogt} which would
finish the proof. However, there is an alternative direct and elementary argument which avoids appealing to
\emph{loc.~cit.} and goes as follows. First notice that, by the same argument as in the proof of
Theorem~\sref{thm:algebraic} point~\eqref{structure-pt:algebraic:sub}, the $k$-rank of a semi-simple subgroup acting minimally
equals the $k$-rank of $G$ (this holds in all cases, not only in the non-Archimedean one). Therefore, the
inclusion of spherical buildings $\mathscr{B}\bar L \to \mathscr{B} G$ provided by the group inclusion $\bar L
\to G$ has the property that $\mathscr{B}\bar L$ is a top-dimensional sub-building of $ \mathscr{B} G$. An
elementary argument (see~\cite[Lemma~3.3]{KleinerLeeb_2006}) shows that the union $Y$ of all apartments of $X$
bounded by a sphere in $\mathscr{B}\bar L$ is a closed convex subset of $X$. Clearly $Y$ is $\bar L$-invariant,
hence $Y = X$ by minimality. Therefore $\mathscr{B}\bar L = \mathscr{B} G$, which finally implies that $\bar L =
G$.
\end{proof}

\subsection{The limit set of subgroups of finite covolume}\label{sec:full:limit}
Let $X$ be a complete \cat space and $G$ a group acting by isometries on $X$. Recall that the \textbf{limit
set}\index{limit set} $\Lambda G$ of $G$ is the intersection of the boundary $\bd X$ with the closure of the
orbit $G.x_0$ in $\overline{X}=X\sqcup \bd X$ of any $x_0\in X$, this set being independent of $x_0$.

\begin{prop}\label{prop:LimitSet}
Let $G$ be a locally compact group acting continuously by isometries on a complete \cat space $X$. If $\Gamma<G$ is any closed subgroup
with finite invariant covolume, then $\Lambda \Gamma = \Lambda G$.
\end{prop}

Consider the following immediate corollary, which in the special case of Hadamard manifolds follows from
the \emph{duality condition}, see~1.9.16 and 1.9.32 in~\cite{Eberlein}.

\begin{cor}\label{cor:LimitSetLattice}
Let $G$ be a locally compact group with a continuous action by isometries on a proper \cat space. If the $G$-action is cocompact,
then any lattice in $G$ has full limit set in $\bd X$.\hfill\qedsymbol
\end{cor}

\begin{proof}[Proof of Proposition~\ref{prop:LimitSet}]
We observe that for any non-empty open set $U\se G$ there is a compact set $C\se G$ such that
$U^{-1} \Gamma C = G$. Indeed, (using an idea of Selberg, compare Lemma~1.4 in~\cite{Borel60}), it suffices to take
$C$ so large that
$$\mu\big(\Gamma C) > \mu(\Gamma\bsl G) - \mu(\Gamma U),$$
where $\mu$ denotes an invariant measure on $\Gamma\bsl G$; any right translate of $\Gamma U$ in $\Gamma\bsl G$ will then meet
$\Gamma C$.

Now let $\xi\in \Lambda G$ and $x_0\in X$. For any neighbourhood $V$ of $\xi$ in $\bd X$, we shall construct an
element in $\Lambda \Gamma\cap V$. Let $U\se G$ be a compact neighbourhood of the identity in $G$ such that
$U\xi\se V$ and let $\{g_n\}$ be a sequence of elements of $G$ with $g_n x_0$ converging to $\xi$ (one uses nets
if $X$ is not separable). In view of the above observation, there are sequences $\{u_n\}$ in $U$ and $\{c_n\}$
in $C$ such that $u_n g_n c_n^{-1}\in\Gamma$. The points $g_n c_n^{-1} x_0$ remain at bounded distance of $g_n
x_0$ as $n\to\infty$, and thus converge to $\xi$. Therefore, choosing an accumulation point $u$ of $\{u_n\}$ in
$U$, we see that $u \xi$ is an accumulation point of $\{u_n g_n c_n^{-1} x_0\}$, which is a sequence in $\Gamma
x_0$.
\end{proof}

For future use, we observe a variant of the above reasoning yielding a more precise fact in a simpler situation:

\begin{lem}\label{lem:compact_ray}
Let $G$ be a locally compact group with a continuous cocompact action by isometries on a proper \cat space $X$.
Let $\Gamma<G$ be a lattice and $c:\RR_+\to X$ a geodesic ray such that $G$ fixes $c(\infty)$. Then there
is a sequence $\{\gamma_i\}$ in $\Gamma$ such that $\gamma_i c(i)$ remains bounded over $i\in \NN$.
\end{lem}

\begin{proof}
For the same reason as above, there is a compact set $U\se G$ such that $G=U\Gamma U^{-1}$.
Choose now $\{g_i\}$ such that $g_i c(i)$ remains bounded and write $g_i = u_i \gamma_i v_i^{-1}$ with $u_i, v_i\in U$.
We have
\begin{multline*}
d(\gamma_i c(i), c(0)) = d(g_i v_i c(i), u_i c(0)) \leq d(g_i v_i c(i), g_i c(i)) + d(g_i c(i),u_i c(0))\\
\leq d(v_i c(i), c(i)) + d(g_i c(i), c(0)) + d(u_ic(0), c(0)).
\end{multline*}
This is bounded independently of $i$ because $d(v_i c(i), c(i)) \leq d(v_i c(0), c(0))$ since  $c(\infty)$
is $G$-fixed.
\end{proof}

We shall also need the following:
\begin{lem}\label{lem:UGammaU}
A locally compact group containing a finitely generated subgroup whose closure has finite covolume is compactly
generated.
\end{lem}
\begin{proof}
Denoting the closure of the given finitely generated subgroup by $\Gamma$, we can write $G = U \Gamma C$ as in
the proof of Proposition~\ref{prop:LimitSet} with both $\overline{U}$ and $C$ compact. Since $\Gamma$ is a locally compact
group containing a finitely generated dense subgroup, it is compactly generated and the conclusion follows.
\end{proof}

\section{\cat lattices, I: the Euclidean factor}
\subsection{Preliminaries on lattices}
We begin this section with a few well known basic facts about general lattices.

\begin{prop}\label{prop:Raghunathan}
Let $G$ be a locally compact second countable group and $N \lhd G$ be a closed normal subgroup.
\begin{enumerate}
\item Given a closed cocompact subgroup $\Gamma < G$, the projection of $\Gamma$ on $G/N$ is closed if and only if
$\Gamma \cap N$ is cocompact in $N$.

\item Given a lattice $\Gamma < G$,  the projection of $\Gamma$ on $G/N$ is discrete if and only if $\Gamma \cap N$
is a lattice in $N$.
\end{enumerate}
\end{prop}

\begin{proof}
See Theorem~1.13 in~\cite{Raghunathan}.
\end{proof}

The second well known result is straightforward to establish:

\begin{lem}\label{lem:CompactOpenTrick}
Let $G= H \times D$ be a locally compact group. Given a lattice $\Gamma < G$
and a compact open subgroup $Q < D$, the subgroup $\Gamma_Q := \Gamma \cap (H \times Q)$ is a lattice in $H
\times Q$, which is commensurated by $\Gamma$.

If moreover $G/\Gamma$ is compact, then so is $(H\times Q)/\Gamma_Q$.\hfill\qed
\end{lem}

\noindent
(As we shall see in Lemma~\ref{lem:CommensuratorTrick} below, there is a form of converse.)

\medskip

Let $X$ be a proper \cat space and $G=\Isom(X)$ be its isometry group. Given
a discrete group $\Gamma$ acting properly and cocompactly on $X$, then the quotient $G\backslash X$ is compact
and the image of $\Gamma$ in $G$ is a cocompact lattice (note that the kernel of the map $\Gamma \to \Isom(X)$
is finite). Conversely, if the quotient $G\backslash X$ is compact, then any cocompact lattice of $G$ is a
discrete group acting properly and cocompactly on $X$.

\begin{lem}\label{lem:finiteGeneration}
In the above setting, $G$ is compactly generated and $\Gamma$ is finitely generated.
\end{lem}

\begin{proof}
For lack of finding a classical reference, we refer to Lemma~22 in~\cite{Mineyev-Monod-Shalom}).
\end{proof}

\subsection{Variations on Auslander's theorem}

\begin{lem}\label{lem:Auslander}
Let $A = \RR^n\rtimes O(n)$ and $S$ be a semi-simple Lie group without compact factor. Any lattice $\Gamma$ in $G
= A \times S$ has a finite index subgroup $\Gamma^0$ which splits as a direct product $\Gamma^0 \cong \Gamma_A
\times \Gamma'$, where $\Gamma_A = \Gamma \cap (A\times 1)$ is a lattice in $(A\times 1)$.
\end{lem}

\begin{proof}
Let $V = \RR^n$ denote the translation subgroup of $A$ and $U$ denote the closure of the projection of $\Gamma$
to $S$. The subgroup $U < S$ is closed of finite covolume; therefore it is either discrete or it contains a
semi-simple subgroup of positive dimension by Borel's density theorem (in fact one could be more precise using
the Main Result of~\cite{Prasad77}, but this is not necessary for the present purposes). On the other hand,
Auslander's theorem~\cite[Theorem~8.24]{Raghunathan} ensures that the identity component of the projection of
$\Gamma$ in $S \times A/V$ is soluble, from which it follows that $U$ has a connected soluble normal subgroup.
Thus $U$ is discrete. Therefore, by Proposition~\ref{prop:Raghunathan}, the group $\Gamma_A = \Gamma \cap
(A\times 1)$ is a lattice in $(A\times 1)$. In particular $\Gamma_A$ is virtually
Abelian~\cite[Corollary~4.1.13]{Thurston}.

Since the projection of $\Gamma$ to $S$ is a lattice in $S$, it is finitely generated~\cite[6.18]{Raghunathan}.
Therefore $\Gamma$ possesses a finitely generated subgroup $\Lambda$ containing $\Gamma_A$ and whose projection
to $S$ coincides with the projection of $\Gamma$. Notice that $\Lambda$ is a lattice in $S \times A$ by
\cite[Theorem~23.9.3]{Simonnet}; therefore $\Lambda$ has finite index in $\Gamma$, which shows that $\Gamma$ is
finitely generated.

Since $\Gamma_A$ is normal in $\Gamma$, the projection $\Gamma^A$ of $\Gamma$ to $A$ normalises the lattice
$\Gamma_A$ and is thus virtually Abelian. Hence $\Gamma^A$ is a finitely generated virtually Abelian group which
normalises $\Gamma_A$. Therefore $\Gamma^A$ has a finite index subgroup which splits as a direct product of the
form $\Gamma_A \times C$, and the preimage $\Gamma'$ of $C$ in $\Gamma$ is a normal subgroup which intersects
$\Gamma_A$ trivially. In particular the group $\Gamma' \cdot \Gamma_A \cong \Gamma' \times \Gamma_A$ is a finite
index normal subgroup of $\Gamma$, as desired.
\end{proof}

\begin{lem}\label{lem:AbstractCommensurator}
Let $\Gamma$  be a group containing a subgroup of the form  $\Gamma^0 \cong \Gamma^0_S \times \Gamma^0_A$, where
$\Gamma^0_S$ is isomorphic to a lattice in a semi-simple Lie group with trivial centre and no compact factor, and
$\Gamma^0_A$ is amenable. If $\Gamma$ commensurates $\Gamma^0$, then $\Gamma$ commensurates both $\Gamma^0_S$
and $\Gamma^0_A$.
\end{lem}
\begin{proof}
Let $\Gamma^1 \cong \Gamma^1_S \times \Gamma^1_A$ be a conjugate of $\Gamma^0$ in $\Gamma$. The projection of
$\Gamma^0 \cap \Gamma^1$ to $\Gamma^0_S$ is a finite index subgroup of $\Gamma^0_S$. By Borel density theorem,
it must therefore have trivial amenable radical. In particular the projection of $\Gamma^0 \cap \Gamma^1_A$ to
$\Gamma^0_S$ is trivial. Therefore the image of the projection of $\Gamma^0 \cap \Gamma^1_S$ (resp. $\Gamma^0
\cap \Gamma^1_A$) to $\Gamma^0_S$ (resp. $\Gamma^0_A$) is of finite index. The desired assertion follows.
\end{proof}

\begin{prop}\label{prop:SxAxD}
Let $A = \RR^n\rtimes O(n)$, $S$ be a semi-simple Lie group with trivial centre and no compact factor, $D$ be a
totally disconnected
locally compact group and $G = S \times A \times D$. 
Then any finitely generated lattice $\Gamma < G$ has a finite index subgroup $\Gamma_0$ which splits as a direct
product $\Gamma_0 \cong \Gamma_A \times \Gamma'$, where $\Gamma_A \se \Gamma \cap (1 \times A \times D)$
is a finitely generated virtually Abelian subgroup whose projection to $A$ is a lattice.
\end{prop}

\begin{proof}
Let $Q < D$ be a compact open subgroup. By Lemma~\ref{lem:CompactOpenTrick}, the intersection $\Gamma^0 = \Gamma
\cap (S \times A \times Q)$ is a lattice in $S \times A \times Q$, which is commensurated by $\Gamma$. Since $Q$
is compact, the projection of $\Gamma^0$ to $S \times A$ is a lattice, to which we may apply
Lemma~\ref{lem:Auslander}. Upon replacing $\Gamma^0$ by a finite index subgroup (which amounts to replacing $Q$
by an open subgroup), this yields two normal subgroups $\Gamma^0_S,  \Gamma^0_A < \Gamma^0$ and a decomposition
$\Gamma^0 = \Gamma^0_S \cdot \Gamma^0_A$, where  $\Gamma^0_S \cap \Gamma^0_A \se Q$ and $\Gamma^0_A = \Gamma^0
\cap (1 \times A \times Q)$ is a finitely generated virtually Abelian group whose projection to $A$ is a
lattice.

By virtue of Lemma~\ref{lem:AbstractCommensurator}, we deduce that the image of the projection of $\Gamma$ to
$A$ commensurates a lattice in $A$. But the commensurator of any lattice in $A$ is virtually Abelian. Therefore,
upon replacing $\Gamma$ by a finite index subgroup, it follows that the projection of $\Gamma$ to $A$ normalises
the projection of $\Gamma^0_A$. We now define
$$
\Gamma_A = \bigcap_{\gamma \in \Gamma} \gamma \Gamma^0_A \gamma \inv.
$$
Then the projection of $\Gamma_A$ coincides with the projection of $\Gamma^0_A$ since $A$ is Abelian; in particular it is still a
lattice. Furthermore, the subgroup $\Gamma_A$ is normal in $\Gamma$. We now proceed as in the proof of
Lemma~\ref{lem:Auslander}. Since the projection of $\Gamma$ to $A$ is finitely generated and virtually Abelian,
we may thus find in this group a virtual complement to the image of the projection of $\Gamma_A$. Let $\Gamma'$
be the preimage of this complement in $\Gamma$. Then, upon replacing $\Gamma$ by a finite index subgroup, the
group $\Gamma'$ is normal in $\Gamma$ and $\Gamma = \Gamma_A\cdot \Gamma'$. Since $\Gamma_A$ is normal as well,
the commutator $[\Gamma_A, \Gamma']$ is contained in the intersection $\Gamma_A \cap \Gamma'$, which is trivial
by construction. This finally shows that $\Gamma \cong \Gamma_A \times \Gamma'$, as desired.
\end{proof}

\begin{remark}\label{rem:SxAxD}
In the setting of Proposition~\ref{prop:SxAxD}, assume that any compact subgroup of $D$ normalised by $\Gamma$ is trivial.
Then $\Gamma_A\se 1\times A\times 1$ and the projection of $\Gamma$ to $S\times D$ is discrete. Indeed, the definition of
$\Gamma_A$ given in the proof shows that it is contained in $1\times A\times \gamma Q \gamma\inv$ for all $\gamma\in\Gamma$
and under the current assumptions the intersection $\bigcap \gamma Q \gamma\inv$ is trivial. The claim about the projection
to $S\times D$ follows from Proposition~\ref{prop:Raghunathan}.
\end{remark}

\subsection{Lattices, the Euclidean factor and fixed points at infinity}
Given a proper \cat space $X$ and a discrete group $\Gamma$ acting properly and cocompactly, it is a well known
open question, going back to M.~Gromov~\cite[\S\,$6.\mathrm{B}_3$]{Gro93}, to determine whether the presence of an
$n$-dimensional flat in $X$ implies the existence of a free Abelian group of rank $n$ in $\Gamma$.
(In the manifold case, see problem~65 on Yau's list~\cite{YauPB}.)
Here we propose the following theorem; the special case where $X/\Gamma$ is a compact Riemannian manifold is the
main result of Eberlein's article~\cite{Eberlein83} (compare also the earlier Theorem~5.2 in~\cite{Eberlein80}).

\begin{thm}\label{thm:lattices:EuclideanSplitting}
Let $X$ be a proper \cat space such that $G = \Isom(X)$ acts cocompactly. Suppose that $X\cong\RR^n\times X'$.

\begin{enumerate}
\item Any finitely generated lattice $\Gamma < G$ has a finite index subgroup $\Gamma_0$ which splits as a direct
product $\Gamma_0 \cong \ZZ^n \times \Gamma'$.\label{pt:lattices:EuclideanSplitting}

\item If moreover $X$ is $G$-minimal (\emph{e.g.} if $X$ is geodesically complete), then $\ZZ^n$ acts trivially
on $X'$ and as a lattice on $\RR^n$; the projection of $\Gamma$ to $\Isom(X')$ is
discrete.\label{pt:lattices:EuclideanSplitting:min}
\end{enumerate}
\end{thm}

\noindent We recall that cocompact lattices are automatically finitely generated in the above setting,
Lemma~\ref{lem:finiteGeneration}. The following example shows that, without the assumption  that $G$ acts
minimally, the projection of $\Gamma$ to $\Isom(X')$ should not be expected to have discrete image:

\begin{example}\label{ex:EuclidienPenché}
Let $X$ be the closed submanifold of $\RR^3$ defined by $X = \{(x, y, z) \in \RR^3 \; | \; 1 \leq z \leq 2\}$
and consider the following Riemannian metric on $X$:
$$ds^2 =  dx^2 + z^2 dy^2 + dz^2.$$
One readily verifies that it is non-positively curved; thus $X$ is a \cat manifold. Clearly $X$ splits off a
one-dimensional Euclidean factor along the $x$-axis. Moreover the group $H \cong \RR^2$ of all translations
along the $xy$-plane preserves $X$ and acts cocompactly. Let $\Gamma$ be the subgroup of $H$ generated by $a$
and $b$, where
$$a : (x, y, z) \mapsto (x, y, z)+(\sqrt{2}, 1, 0) \hspace{.5cm} \text{and} \hspace{.5cm}%
b : (x, y, z) \mapsto (x, y, z)+(1, \sqrt{2}, 0).
$$
Then $\Gamma \cong \ZZ^2$ is a cocompact lattice in $\Isom(X)$, but no non-trivial subgroup of $\Gamma$ acts
trivially on the $yz$-factor of $X$. The projection of $\Gamma$ to the isometry group of that factor is not
discrete (see Proposition~\ref{prop:Raghunathan}(ii)).
\end{example}

The above result is the converse to the Flat Torus Theorem when it is stated as
in~\cite[II.7.1]{Bridson-Haefliger}. In particular we deduce that the dimension of the Euclidean de Rham factor
is an invariant of $\Gamma$. In the manifold case, again, this is the main point of~\cite{Eberlein83}.

\begin{cor}\label{cor:lattices:EuclideanSplitting}
Let $X$ be a proper \cat space such that $G = \Isom(X)$ acts cocompactly and minimally.
Let $\Gamma<G$ be a finitely generated lattice.

Then the dimension of the Euclidean factor of $X$ equals the maximal rank of a free Abelian normal subgroup of
$\Gamma$.
\end{cor}

In order to apply Theorem~\sref{thm:Decomposition} and Addendum~\sref{addendum} towards
Theorem~\ref{thm:lattices:EuclideanSplitting}, we will need the following.

\begin{thm}\label{thm:Lattice=>NoFixedPointAtInfty}
Let $X$ be a proper \cat space such that $G= \Isom(X)$ acts cocompactly and contains a finitely generated lattice.
Then $X$ contains a canonical closed convex $G$-invariant $G$-minimal subset $X'\neq \varnothing$
which has no $\Isom(X')$-fixed point at infinity.
\end{thm}

Consider the immediate corollary.

\begin{cor}\label{cor:Lattice=>NoFixedPointAtInfty}
Let $X$ be a proper \cat space such that $G= \Isom(X)$ acts cocompactly and minimally. If $G$ contains a
finitely generated lattice, then $G$ has no fixed point at infinity.\qed
\end{cor}

This shows that the mere existence of a finitely generated lattice imposes restrictions on cocompact \cat
spaces; much more detailed results in that spirit will be given in Section~\ref{sec:geometric:arith}.

We do not know whether the statement of Corollary~\ref{cor:Lattice=>NoFixedPointAtInfty} remains true without
the finite generation assumption on the lattice (see Problem~\ref{pb:FiGen} below).

\begin{example}\label{ex:CocompactPtFixeAuBord}
We emphasise that the full isometry group of a cocompact proper \cat space may have global fixed points at
infinity; in fact, the space might even be homogeneous, as it is the case for E.~Heintze's
manifolds~\cite{Heintze74} mentioned earlier. An even simpler way to construct cocompact proper \cat space with
this property is to mimic Example~\sref{ex:Arbre-triangle}: Start from a regular tree $T$, assuming for definiteness that
the valency is three. Replace every vertex by a congruent copy of an isosceles triangle that is not equilateral, in such a
way that its distinguished vertex always points to a fixed point at infinity (of the initial tree). Then the stabiliser $H$
in $\Isom(T)$ of that point at infinity still acts faithfully and cocompactly on the modified space $T'$; the
construction is so that the isometry group of $T'$ is in fact reduced to $H$.
\end{example}

We shall also establish a strengthening of Corollary~\ref{cor:Lattice=>NoFixedPointAtInfty}, which can be viewed
as a form of Borel (or geometric) density theorem without assumption about fixed points at infinity.\index{Borel density}

\begin{thm}\label{thm:Lattice=>NoFixedPointAtInfty_bis}
Let $X$ be a proper \cat space such that $G= \Isom(X)$ acts cocompactly and minimally.
Assume there is a finitely generated lattice $\Gamma< G$. Then $\Gamma$ acts minimally
on $X$ and moreover all $\Gamma$-fixed points at infinity are contained in the
boundary of the (possibly trivial) Euclidean factor of $X$.
\end{thm}

We now turn to the proofs.
In the case of a discrete cocompact group $\Gamma=G$, a version of the following was first established by
Burger--Schroeder~\cite{BurgerSchroeder87} (as pointed out in~\cite[Corollary~2.7]{AB98}).

\begin{prop}\label{prop:LatticeAndFixedPoints}
Let $X$ be a proper \cat space, $G<\Isom(X)$ a closed subgroup whose action on $X$ is cocompact
and $\Gamma<G$ a finitely generated lattice.
Then there exists a $\Gamma$-invariant closed convex subset $Y\se X$ which splits
$\Gamma$-equivariantly as $Y = E \times W$, where $E$ is a (possibly $0$-dimensional) Euclidean space
on which $\Gamma$ acts by translations and such that $\bd E$ contains the fixed point set of $G$ in $\bd X$.
\end{prop}

\begin{proof}
We can assume that there are $G$-fixed points at infinity,
since otherwise there is nothing to prove. We claim that for any $G$-fixed point $\xi$
there is a geodesic line $\sigma:\RR\to X$ with $\sigma(+\infty)=\xi$ such that any
$\gamma\in \Gamma$ moves $\sigma$ to within a bounded distance of
itself~--- and hence to a parallel line by convexity of the metric.

Indeed, let $c:\RR_+\to X$ be a geodesic ray with $c(\infty)=\xi$ and let $\{\gamma_i\}$ be as in
Lemma~\ref{lem:compact_ray}. Then, by Arzel\`a--Ascoli, there is a subsequence $I\se \NN$
and a geodesic line $\sigma:\RR\to X$ such that $\sigma(t)=\lim_{i\in I}\,\gamma_i\, c(t+i)$ for all $t$.
Since each $g\in G$ has bounded displacement along $c$, the sequence $\{\gamma_i g \gamma_i^{-1}\}_{i\in I}$
is bounded and thus we can assume that it converges for all $g$ (recalling that $G$ is second countable,
but we shall only consider $g\in\Gamma$ anyway).
Since $\Gamma$ is discrete and finitely generated, we can further restrict $I$ so that there is
$\gamma_\infty\in\Gamma$ such that
$$\gamma_i \gamma \gamma_i^{-1} = \gamma_\infty \gamma \gamma_\infty^{-1}\kern10mm \forall\,\gamma\in\Gamma, i\in I.$$
Since
\begin{multline*}
d(\gamma\gamma_\infty^{-1}\sigma(t), \gamma_\infty^{-1}\sigma(t)) =
\lim_{i\in I} d(\gamma_i^{-1} \gamma_\infty \gamma \gamma_\infty^{-1}\gamma_i c(t+i), c(t+i))\\
=\lim_{i\in I} d(\gamma c(t+i), c(t+i)) \leq d(\gamma c(0), c(0)),
\end{multline*}
It now follows that every $\gamma\in\Gamma$ has bounded displacement length along the geodesic
$\gamma_\infty^{-1}\sigma$. Thus the same holds for the geodesic $\sigma$ which is therefore
(by convexity) translated to a parallel line by each element of $\Gamma$ as claimed.

Consider a flat $E\se X$ that is maximal for the property that each element of $\Gamma$ has constant
displacement length on $E$. Let $Y$ be the union of all flats that are at finite distance from $E$. One shows
that $Y$ splits as $Y\cong E\times W$ for some closed convex $W\se X$ using the Sandwich
Lemma~\cite[II.2.12]{Bridson-Haefliger} and Lemma~II.2.15 of~\cite{Bridson-Haefliger} just like in
Section~\sref{sec:largeRadius}. The definition of $Y$ shows that $\Gamma$ preserves $Y$ as well as its
splitting and acts on the $E$ coordinate by translations.

It remains to show that any $G$-fixed point $\xi\in\bd X$ belongs to $\bd E$. First, $\xi\in\bd Y$ since $\bd
Y=\bd X$ by Corollary~\ref{cor:LimitSetLattice}; we thus represent $\xi$ by a ray $c:\RR_+\to Y$. Let now
$\sigma$ be a geodesic line as provided by the claim. We can assume that $\sigma$ lies in $Y$ because it was
constructed from $\Gamma$-translates of $c$ and $Y$ is $\Gamma$-invariant. One can write $\sigma=(\sigma_E,
\sigma_W)$ where $\sigma_E$, $\sigma_W$ are linearly re-parametrised geodesics in $E$ and $W$,
see~\cite[I.5.3]{Bridson-Haefliger}. We need to prove that $\sigma_W$ has zero speed. Since any given $\gamma\in
\Gamma$ has constant displacement along $\sigma$ and on each of the parallel copies of $E$ individually, its
displacement is constant on the union of all parallel copies of $E$ visited by $\sigma$, which is $E\times
\sigma_W(\RR)$. The latter being again a flat, the maximality of $E$ shows that $\sigma_W$ is constant.
\end{proof}

\begin{proof}[Proof of Theorem~\ref{thm:Lattice=>NoFixedPointAtInfty}]
Let $Y=E\times W\se X$ be as in Proposition~\ref{prop:LatticeAndFixedPoints}. Recall that $\bd Y=\bd X$ by
Corollary~\ref{cor:LimitSetLattice}. We claim that $\bd X$ has circumradius~$>\pi/2$. Indeed, it would otherwise
have a $G$-fixed circumcentre by Proposition~\sref{prop:BalserLytchak}, but this circumcentre cannot belong to
$\bd E$ since $E$ is Euclidean; this contradicts Proposition~\ref{prop:LatticeAndFixedPoints}. We now apply
Corollary~\sref{cor:CanonicalFullSubset}. This yields a canonical $G$-invariant closed convex subset $X'$, which
is minimal with respect to the property that $\bd X' = \bd X$. It follows in particular by
Corollary~\ref{cor:LimitSetLattice} that $\Gamma$ acts minimally on $X'$. Let now $X' = E' \times X'_0$ be the
canonical splitting, where $E'$ is the maximal Euclidean factor~\cite[II.6.15]{Bridson-Haefliger}. On the one
hand, since $X'$ is $\Gamma$-minimal, Proposition~\ref{prop:LatticeAndFixedPoints} applied to $X'$ shows that
$G$ has no fixed points in $\bd X'_0$ since $E'$ is maximal as a Euclidean factor. On the other hand,
$\Isom(E')$ fixes no point at infinity on $E'$. We deduce that $\Isom(X') \cong \Isom(E') \times \Isom(X'_0)$
has indeed no fixed point at infinity.
\end{proof}

\begin{proof}[End of proof of Theorem~\ref{thm:Lattice=>NoFixedPointAtInfty_bis}]
Arguing as in the proof of Theorem~\ref{thm:Lattice=>NoFixedPointAtInfty}, we establish that $X$ is $\Gamma$-minimal.
Let $X=X'\times E$ be the canonical splitting, where $E$ is the maximal Euclidean factor.
Since any isometry of $X$ decomposes uniquely as isometries of $E$ and $X'$ (II.6.15 in~\cite{Bridson-Haefliger}),
is suffices to show that $\Gamma$ has no fixed point in $\bd X'$. This follows from Proposition~\ref{prop:coAmenable}
applied to the $G$-action on $X'$.
\end{proof}

\begin{proof}[End of proof of Theorem~\ref{thm:lattices:EuclideanSplitting}]
Assume first that $X$ is $G$-minimal, recalling that this is the case if $X$ is geodesically complete by
Lemma~\sref{lem:cocompact:minimal}. In view of Corollary~\ref{cor:Lattice=>NoFixedPointAtInfty}, we can apply
Theorem~\sref{thm:Decomposition} and we are therefore in the setting of Proposition~\ref{prop:SxAxD}.
Since the group $\Gamma_A$ provided by that proposition contains a finite index subgroup isomorphic to $\ZZ^n$,
we have already established~\eqref{pt:lattices:EuclideanSplitting} under the additional minimality assumption.

\smallskip

In order to show~\eqref{pt:lattices:EuclideanSplitting:min}, it suffices by Remark~\ref{rem:SxAxD} to prove that any
compact subgroup of $G$ normalised by $\Gamma$ is trivial. This follows from the fact that $X$ is $\Gamma$-minimal,
as established in Theorem~\ref{thm:Lattice=>NoFixedPointAtInfty_bis}.

\smallskip

It remains to prove~\eqref{pt:lattices:EuclideanSplitting} without the assumption that $X$ is $G$-minimal.
Let $Y\se X$ be the $G$-minimal set provided by Theorem~\ref{thm:Lattice=>NoFixedPointAtInfty} and let
$Y\cong \RR^m\times Y'$ be its Euclidean decomposition. Then $m\geq n$ because of the characterisation of
the Euclidean factor in terms of Clifford isometries~\cite[II.6.15]{Bridson-Haefliger}; indeed, any (non-trivial)
Clifford isometry of $X$ restricts non-trivially to $Y$ because $Y$ has finite co-diameter. The kernel $F\lhd\Gamma$
of the $\Gamma$-action on $Y$ is finite and thus we can assume that it is central upon replacing $\Gamma$ with a finite
index subgroup. Passing to a further finite index subgroup, we know from the minimal case that $\Gamma/F$ splits as
$\Gamma/F = \ZZ^m \times \Lambda'$. Let $\Gamma_{\ZZ^m}, \Gamma'\lhd \Gamma$ be the pre-images in $\Gamma$ of those factors.
Thus we can write $\Gamma=\Gamma_{\ZZ^m}\cdot\Gamma'$ with $\Gamma_{\ZZ^m}\cap\Gamma'\se F$. It is straightforward that a finite
central extension of $\ZZ^m$ is virtually $\ZZ^m$ (see \emph{e.g.}~\cite[II.7.9]{Bridson-Haefliger}). Therefore $\Gamma$
contains a finite index subgroup isomorphic to $\ZZ^m\times \Gamma'$ and the result follows since $m\geq n$.
\end{proof}

\begin{proof}[Proof of Corollary~\ref{cor:lattices:EuclideanSplitting}]
Notice that a splitting $\Gamma_0\cong \ZZ^n\times \Gamma'$
with $\Gamma_0$ normal and $n$ maximal provides a normal subgroup $\ZZ^n\lhd \Gamma$ since $\ZZ^n$ is characteristic in
$\Gamma_0$. Therefore, given Theorem~\ref{thm:lattices:EuclideanSplitting}, it only remains to see that a normal
$\ZZ^n\lhd \Gamma$ of maximal rank forces $X$ to have a Euclidean factor of dimension at least $n$. Otherwise, the projection
of $\Gamma$ to the non-Euclidean factor $X'$ would be a lattice by
Theorem~\ref{thm:lattices:EuclideanSplitting}\eqref{pt:lattices:EuclideanSplitting:min} and contain an infinite normal
amenable subgroup, contradicting Corollary~\ref{cor:LatticeNormaliser}\eqref{pt:LatticeNormaliser:ramen}.
\end{proof}

Finally, we record that Theorem~\ref{thm:fglattices} is contained in
Theorem~\ref{thm:lattices:EuclideanSplitting} and Corollary~\ref{cor:lattices:EuclideanSplitting} for~(i), and
Corollary~\ref{cor:Lattice=>NoFixedPointAtInfty} and Theorem~\ref{thm:Lattice=>NoFixedPointAtInfty_bis}
for~(ii).

\section{\cat lattices, II: products}
\subsection{Irreducible lattices in \cat spaces}
Recall from that a (topological) group is called \textbf{irreducible}\index{irreducible!group}
if no (open) finite index subgroup splits non-trivially as a direct product of (closed) subgroups. For example,
any locally compact group acting continuously, properly, minimally, without fixed point at infinity on an
irreducible proper \cat space is irreducible by Theorem~\sref{thm:geometric_simplicity}.

In particular, an abstract group $\Gamma$ is irreducible if it does not virtually split. This terminology is
inspired by the concept of irreducibility for closed manifolds, which means that no finite cover of the manifold
splits non-trivially. Of course, the universal cover of such a manifold can still split. Indeed, one gets many
classical \cat groups by considering ``irreducible lattices'' in products of simple Lie groups or more
generally of semi-simple algebraic groups over various local fields.

\smallskip

The latter concept of irreducibility for lattices is defined as follows: A lattice $\Gamma<G=G_1 \times \dots
\times G_n$ in a product of locally compact groups is called an \textbf{irreducible
lattice}\label{def:IrreducibleLattice}\index{irreducible!lattice} if its projections to any subproduct of the
$G_i$'s are dense and each $G_i$ is non-discrete.

\smallskip

The point of this notion (and of the nearly confusing terminology) is that it prevents $\Gamma$
and its finite index subgroups from splitting as a product of lattices in $G_i$. Moreover, if all
$G_i$'s are centre-free simple Lie (or algebraic) groups without compact factors, the irreducibility
of $\Gamma$ as a lattice is equivalent to its irreducibility as a group in and for itself;
this is a result of Margulis~\cite[II.6.7]{Margulis}. As we shall see in Theorem~\ref{thm:lattices:Irred}
below, a version of this equivalence holds for lattices in the isometry group of a \cat space.

\begin{remark}\
\begin{enumerate}
\item The non-discreteness of $G_i$ is often omitted from this definition; the difference is inessential since
the notion of a lattice is trivial for discrete groups. Notice however that our definition ensures that all
$G_i$ are non-compact and that $n\geq 2$.

\item One verifies that \emph{any} lattice $\Gamma<G=G_1\times G_2$ is an irreducible lattice in the product $G^*<G$ of the
closures $G^*_i<G_i$ of its projections to $G_i$ (provided these projections are non-discrete).
\end{enumerate}
\end{remark}

The following geometric version of Margulis' criterion contains Theorem~\ref{thm:irred-intro}
from the Introduction.

\begin{thm}\label{thm:lattices:Irred}
Let $X$ be a proper \cat space, $G< \Isom(X)$ a closed subgroup acting cocompactly  on $X$, and $\Gamma < G$ a
finitely generated lattice.

\begin{itemize}
\item[(i)] If $\Gamma$ is irreducible as an abstract group, then for any finite index subgroup $\Gamma_0<\Gamma$
and any $\Gamma_0$-equivariant splitting $X = X_1 \times X_2$ with $X_1$ and $X_2$ non-compact, the projection
of $\Gamma_0$ to both $\Isom(X_i)$ is non-discrete.

\item[(ii)] If in addition the $G$-action is minimal, then the converse statement holds as well.
\end{itemize}
\end{thm}

\begin{remark}
Recall that the $G$-minimality is automatic if $X$ is geodesically complete (Lemma~\sref{lem:cocompact:minimal}).
Statement~(ii) fails completely without minimality (as witnessed for instance by the uncosmopolitan mien of an
equivariant mane).
\end{remark}

\begin{proof}[Proof of Theorem~\ref{thm:lattices:Irred}]
Suppose $\Gamma$ irreducible. Let $X' \se X$ be the canonical subspace provided by
Theorem~\ref{thm:Lattice=>NoFixedPointAtInfty}. By Theorem~\ref{thm:lattices:EuclideanSplitting}, the space $X'$
has no Euclidean factor unless $X=\RR$ and $\Gamma$ is virtually cyclic, in which case the desired statement is
empty.

We first deal with the case when $G$ acts minimally on $X$; by Theorem~\ref{thm:density} this amounts to assume
$X = X'$. Suppose for a contradiction that for $\Gamma_0$ and $X' = X'_1 \times X'_2$ as in the statement, the
projection $G_1$ of $\Gamma_0$ to $\Isom(X'_1)$ is discrete. Let $G_2$ be the closure of the projection of
$\Gamma_0$ to $\Isom(X'_2)$ and notice that both $G_i$ are compactly generated since $\Gamma$ and hence also
$\Gamma_0$ is finitely generated. The projection $\Gamma_2$ of $\Gamma_0\cap(1\times G_2)$ to $G_2$ is a lattice
(by Lemma~\ref{lem:CompactOpenTrick} or by Proposition~\ref{prop:Raghunathan}); being normal, it is cocompact
and hence finitely generated. By Theorem~\ref{thm:Lattice=>NoFixedPointAtInfty_bis}, the group $\Gamma_0$, and
hence also $G_2$, acts minimally and without fixed point at infinity on $X'_2$. Therefore
Corollary~\ref{cor:LatticeNormaliser}\eqref{pt:LatticeNormaliser:centra} implies that the centraliser
$\centra_{G_2}(\Gamma_2)$ is trivial. But $\Gamma_2$ is discrete, normal in $G_2$, and finitely generated. Hence
$\centra_{G_2}(\Gamma_2)$ is open and thus $G_2$ is discrete. Therefore, the product $G_1\times G_2$, which
contains $\Gamma_0$, is a lattice in $\Isom(X'_1) \times \Isom(X'_2)$ and thus in $G$. Now the index of
$\Gamma_0$ in $G_1 \times G_2$ is finite and thus $\Gamma_0$ splits virtually, a contradiction.

We now come back to the general case $X' \se X$ and suppose that $X$ possesses a $\Gamma_0$-equivariant
splitting $X = X_1 \times X_2$. The group $H = \Isom(X_1) \times \Isom(X_2) < \Isom(X)$ contains $\Gamma_0$;
hence its action on $X'$ is minimal without fixed point at infinity by
Corollary~\ref{cor:Lattice=>NoFixedPointAtInfty}. Therefore, the splitting
theorem~\cite[Theorem~9]{Monod_superrigid} implies that $X'$ possesses a $\Gamma_0$-equivariant splitting $X' =
X'_1 \times X'_2$ induced by $X = X_1 \times X_2$ via $H$. Upon replacing $\Gamma_0$ be a finite index subgroup,
the preceding paragraph thus yields a  splitting $\Gamma_0/F \cong G_1 \times G_2$ of the image of $\Gamma_0$ in
$\Isom(X')$, where $F$ denotes the kernel of the $\Gamma_0$-action on $X'$. Since $F$ is finite, so is the
projection to $\Isom(X_{3-i})$ of the preimage $\widehat{G}_i$ of $G_i$ in $\Gamma$, for $i =1, 2$. Therefore
upon passing to a finite index subgroup we may and shall assume that $\widehat{G}_i$ acts trivially on
$\Isom(X_{3-i})$. Now the subgroup of $\Isom(X_1)\times \Isom(X_2)$ generated by $\widehat{G}_1$ and
$\widehat{G}_2$ splits as $\widehat{G}_1 \times \widehat{G}_2$ and is commensurable to $\Gamma_0$, a
contradiction.

\medskip

Conversely, suppose now that the $G$-action is minimal and that $\Gamma= \Gamma'\times \Gamma''$ splits
non-trivially (after possibly having replaced it by a finite index subgroup). If $X=\RR^n$, then reducibility of
$\Gamma$ forces $n\geq 2$ and we are done in view of the structure of Bieberbach groups. If $X$ is not Euclidean
but has a Euclidean factor, then
Theorem~\ref{thm:lattices:EuclideanSplitting}\eqref{pt:lattices:EuclideanSplitting:min} provides a discrete
projection of $\Gamma$ to the non-Euclidean factor $\Isom(X')$; furthermore, $X'$ is indeed non-compact as
desired since otherwise by minimality it is reduced to a point, contrary to our assumption.

If on the other hand $X$ has no Euclidean factor, then $\Gamma$ acts minimally and without fixed point at
infinity by Theorem~\ref{thm:Lattice=>NoFixedPointAtInfty}. Then the desired splitting is provided by the
splitting theorem~\cite[Theorem~9]{Monod_superrigid}. Both projections of $\Gamma$ are discrete, indeed isomorphic
to $\Gamma'$ respectively $\Gamma''$ because the cited splitting theorem ensures componentwise action of
$\Gamma$.
\end{proof}

We now briefly turn to uniquely geodesic spaces and to the analogues in this setting of some of P.~Eberlein's results
for Hadamard manifolds.

\begin{thm}\label{thm:nonuniform:unique_geod}
Let $X$ be a proper \cat space with uniquely extensible geodesics such that $\Isom(X)$ acts
cocompactly on $X$.

If $\Isom(X)$ admits a finitely generated non-uniform irreducible lattice, then $X$ is a symmetric space
(without Euclidean factor).
\end{thm}

\begin{proof}
The action of $\Isom(X)$ is minimal by Lemma~\sref{lem:cocompact:minimal} and without fixed point at infinity
by Corollary~\ref{cor:Lattice=>NoFixedPointAtInfty}. Thus we can apply Theorem~\sref{thm:Decomposition} and
Addendum~\sref{addendum}. Notice that $\Isom(X)$ itself is non-discrete since it contains a non-uniform lattice;
moreover, if it admits more than one factor in the decomposition of Theorem~\sref{thm:Decomposition},
then the latter are all non-discrete by Theorem~\ref{thm:lattices:Irred}. Therefore, we can apply
Theorem~\sref{thm:nondiscrete:nonbranching} to all factors of $X$. It remains only to justify that $X$ has no
Euclidean factor; otherwise, Auslander's theorem (compare also Theorem~\ref{thm:lattices:EuclideanSplitting})
implies $X=\RR$, which is incompatible with the fact that $\Gamma$ is non-uniform.
\end{proof}

The following related result is due to P.~Eberlein in the manifold case (Proposition~4.5 in~\cite{Eberlein82}).
We shall establish another result of the same vein later without assuming that
geodesics are uniquely extensible (see Theorem~\ref{thm:arith:geometric:lin}).

\begin{thm}\label{thm:Eberlein:unique_geod}
Let $X$ be a proper \cat space with uniquely extensible geodesics and $\Gamma < \Isom(X)$ be a discrete
cocompact group of isometries. If $\Gamma$ is irreducible (as an abstract group) and $X$ is reducible, then $X$
is a symmetric space (without Euclidean factor).
\end{thm}

\begin{proof}
One follows line-by-line the proof of Theorem~\ref{thm:nonuniform:unique_geod}. The only difference is that, in
the present context, the non-discreteness of the isometry group of each irreducible factor of $X$ follows from
Theorem~\ref{thm:lattices:Irred} since $X$ is assumed reducible.
\end{proof}

We can now conclude the proof of Theorem~\ref{thm:Hadamard} from the Introduction. The first statement was
established in Theorem~\sref{thm:nondiscrete:nonbranching}. The second follows from
Theorem~\ref{thm:nonuniform:unique_geod} and the third from Theorem~\ref{thm:Eberlein:unique_geod} in the
uniform case, and from Theorem~\ref{thm:nonuniform:unique_geod} in the non-uniform one.\qed

\subsection{The hull of a lattice}\label{sec:hull}
Let $X$ be a proper \cat space $X$ such that $\Isom(X)$ acts cocompactly on $X$. Let $\Gamma <  \Isom(X)$ be a
finitely generated lattice; note that the condition of finite generation is redundant if $\Gamma$ is cocompact by
Lemma~\ref{lem:finiteGeneration}. Theorem~\ref{thm:Lattice=>NoFixedPointAtInfty} provides a canonical
$\Isom(X)$-invariant subspace $X' \se X$ such that $G = \Isom(X')$ has no fixed point at infinity.

\medskip

In this section we shall define the \textbf{hull}\index{hull} $H_\Gamma$ of the lattice $\Gamma$; this is a
locally compact group $H_\Gamma<\Isom(X')$ canonically attached to the situation and containing the image of
$\Gamma$ in $\Isom(X')$.

\bigskip

For simplicity, we first treat the special case where $\Isom(X)$ acts minimally; thus $X'=X$ and $G=\Isom(X)$.
Applying Theorem~\sref{thm:Decomposition} and Addendum~\sref{addendum}, we see in particular that $\Gamma$
possesses a canonical finite index normal subgroup $\Gamma^* = \Gamma \cap G^*$ which is the kernel of the
$\Gamma$-action by permutation on the set of factors in the decomposition given by Addendum~\sref{addendum}.

In the classical case when $X$ is a symmetric space, the closure of the projection of $\Gamma$ to the
isometry group $\Isom(X_i)$ of each factor is an open subgroup \emph{of finite
index}, as soon as $X$ is reducible. This is no longer true in general, even in the case of Euclidean
buildings. In fact, the same $\Gamma$ may (and generally does) occur as lattice in increasingly
large ambient groups $\Gamma<G<G'< G'' < \cdots$. In order to address this issue, we define the hull
as follows. Consider the closed subgroup $H_{\Gamma^*} < G$ which is the direct product of the closure of the
images of $\Gamma^*$ in each of the factors in the decomposition of
Theorem~\sref{thm:Decomposition}. Then set $H_\Gamma = \Gamma \cdot H_{\Gamma^*}$. In other words, we have inclusions
$$\Gamma\ <\ H_\Gamma\ <\ G.$$

The closed subgroup $H_{\Gamma^*}$ is nothing but the hull of the lattice $\Gamma^*$. It coincides with
$H_\Gamma^* = H_\Gamma \cap G^*$. In particular  $H^*_\Gamma = H_{\Gamma^*}$ is a direct product of irreducible
groups satisfying all the restrictions of Theorem~\sref{thm:geometric_simplicity} (except for the possible
Euclidean motion factor), and the image of $\Gamma^*$ in each of these factors is dense.

\begin{remark}\label{rem:hull}
Notice that $\Gamma$ is always a lattice in $H_\Gamma$ (by~\cite[Lemma~1.6]{Raghunathan}). We emphasise that
$H_\Gamma$ is non-discrete and that $\Gamma^*$ is an irreducible lattice in $H_{\Gamma^*}$ (in the sense of
\S\,\ref{def:IrreducibleLattice}) as soon as $\Gamma$ is irreducible as a group and $X$ is reducible; this follows
from Theorem~\ref{thm:lattices:Irred}.
\end{remark}

We now define the hull $H_\Gamma < G$ in the general situation $G=\Isom(X')$ with $X'\se X$ given by
Theorem~\ref{thm:Lattice=>NoFixedPointAtInfty}. Since $\Isom(X) \backslash X$ is cocompact, it follows that $X'$
is $r$-dense in $X$ for some $r > 0$ and the canonical map $\Isom(X) \to G$ is proper. Let $F_\Gamma\lhd\Gamma$
be the finite kernel of the induced map $\Gamma\to G$ and write $\Gamma':=\Gamma/F_\Gamma$. Then \textbf{the
hull of $\Gamma$}\index{hull!of a lattice} is defined by $H_\Gamma:=H_{\Gamma'}$ (reducing to the above case).

In other words, $\Gamma$ sits in $H_\Gamma$ only modulo the canonical finite kernel $F_\Gamma$. In fact,
$F_\Gamma$ is even canonically attached to $\Gamma$ viewed as an abstract group.

\begin{lem}\label{lem:finiteRad}
$F_\Gamma$ is a (necessarily unique) maximal finite normal subgroup of $\Gamma$.
Moreover, $X'$ is $\Gamma'$-minimal.
\end{lem}

\begin{proof}
The $\Gamma'$-action on $X'$ is minimal by an application of Theorem~\ref{thm:Lattice=>NoFixedPointAtInfty_bis}
and therefore every finite normal subgroup of $\Gamma'$ is trivial.
Since moreover the $\Gamma$-action on $X'$ is proper, it follows that a normal subgroup of $\Gamma$ is finite if
and only if it lies in $F_\Gamma$.
\end{proof}

For later references, we record the following expected fact.

\begin{lem}\label{lem:hull}
Assume that $\Gamma$ is irreducible. If $X'$ is reducible, then $H_\Gamma$ contains the identity component of
$G:=\Isom(X')$. In fact $(H_\Gamma)^\circ = G^\circ$ is a semi-simple Lie group with trivial centre and no compact factor.
\end{lem}

\begin{proof}
By Theorem~\ref{thm:lattices:EuclideanSplitting}, the hypotheses on $\Gamma$ imply that $X'$ has no Euclidean
factor. Thus each almost connected factor of $G^*$ is a simple Lie group with trivial centre and no compact
factor. The projection of $\Gamma^*$ to each of these factors is non-discrete by
Theorem~\ref{thm:lattices:Irred} and the assumption made on $X'$. Its closure is semi-simple and Zariski dense
by Theorem~\ref{thm:density} and Proposition~\ref{prop:Karpelevich-Mostow}. The result follows.
\end{proof}
\subsection{On the canonical discrete kernel}

Let $G=G_1\times G_2$ be a locally compact group and $\Gamma < G$ be an irreducible lattice. It follows from
irreducibility that the projection to $G_i$ of the kernel of the projection $\Gamma\to G_{j\neq i}$ is a normal
subgroup of $G_i$. In other words, we have a canonical discrete normal subgroup $\Gamma_i\lhd G_i$ defined by
$$\Gamma_1\ =\ \mathrm{Proj}_{G_1}\big(\Gamma\cap (G_1\times 1)\big)$$
(and likewise for $\Gamma_2$) which we call the \textbf{canonical discrete kernel}\index{canonical!discrete
kernel} of $G_i$ (depending on $\Gamma$). We observe that the image
$$\ud\Gamma\ = \Gamma/(\Gamma_1\cdot\Gamma_2)$$
of $\Gamma$ in the \textbf{canonical quotient}\index{canonical!quotient} $G_1/\Gamma_1\times G_2/\Gamma_2$ is
still an irreducible lattice (see Proposition~\ref{prop:Raghunathan}(ii)) and has the additional property that
it projects injectively into both factors.

In this subsection, we collect some basic facts on lattices in (products of) totally disconnected locally
compact groups, adapting ideas of M.~Burger and Sh.~Mozes (see Propositions~2.1 and~2.2
in~\cite{Burger-Mozes2}).

\begin{prop}\label{prop:ResiduallyFinite}
Let $\Gamma < G = G_1 \times G_2$ be an irreducible lattice. Assume that $G_2$ is totally disconnected,
compactly generated and without non-trivial compact normal subgroup. If $\Gamma$ is residually finite, then
canonical the discrete kernel $\Gamma_2 = \Gamma \cap (1 \times G_2)$ commutes with the discrete residual
$G_2^{(\infty)}$.
\end{prop}

Recall that the \textbf{discrete residual}\index{discrete residual|see{residual}}\index{residual!discrete}
$G^{(\infty)}$ of a topological group $G$ is by definition the intersection of all open normal subgroups.
It is important to remark that, by Corollary~\sref{cor:ResiduallyDiscrete} the discrete residual of a non-discrete
compactly generated locally compact group without non-trivial compact normal subgroup is necessarily
non-trivial.

\begin{proof}[Proof of Proposition~\ref{prop:ResiduallyFinite}]
By a slight abuse of notation, we shall identify $G_2$ with the subgroup $1 \times G_2$ of $G$. Given a finite
index normal subgroup $\Gamma_0 \lhd \Gamma$, the intersection $\Gamma_0 \cap \Gamma_2$ is a discrete normal
subgroup of $G_2$ (by irreducibility), contained as a finite index subgroup in $\Gamma_2$. Thus $G_2$ acts by
conjugation on the finite quotient $\Gamma_2/ \Gamma_0 \cap \Gamma_2$. In particular the kernel of this action
is a finite index closed normal subgroup, which is thus open. Therefore, the discrete residual $G_2^{(\infty)}$
acts trivially on $\Gamma_2/ \Gamma_0 \cap \Gamma_2$. In other words, this means that $[\Gamma_2 ,
G_2^{(\infty)}] \se \Gamma_0 \cap \Gamma_2$.

Assume now that $\Gamma$ is residually finite. The preceding argument then shows that the commutator $[\Gamma_2
, G_2^{(\infty)}]$ is trivial, as desired.
\end{proof}

\begin{prop}\label{prop:DiscreteNormal}
Let $\Gamma < G = G_1 \times G_2$ be a cocompact lattice in a product of compactly generated locally compact
groups. Assume that $G_2$ is totally disconnected  and that the centraliser in $G_1$ of any cocompact lattice of
$G_1$ is trivial. If the discrete kernel $\Gamma_2 = \Gamma \cap (1 \times G_2)$ is trivial, then the
quasi-centre $\QZ(G_2)$ is topologically locally finite.
\end{prop}

\begin{proof}
Let $S \se \QZ(G_2)$ be a finite subset of the quasi-centre. Then $G_2$ possesses a compact open subgroup $U$
which centralises $S$. By Lemma~\ref{lem:CompactOpenTrick} the group $\Gamma_U = \Gamma \cap (G_1 \times U)$ is
a cocompact lattice in $G_1 \times U$. In particular, there is a finite generating set $T \se \Gamma_U$. By a
lemma of Selberg~\cite{Selberg60}, the group $\centra_\Gamma(T)$ is a cocompact lattice in $\centra_G(T)$. But
$\centra_G(T) = \centra_G(\Gamma_U) \se 1 \times G_2$ since the projection of $\Gamma_U$ to $G_1$ is a cocompact
lattice. Since the discrete kernel $\Gamma \cap (1 \times G_2)$ is trivial by hypothesis, the centraliser
$\centra_\Gamma(T)$ is trivial and, hence, $\centra_G(T)$ is compact. By construction $S$ is contained in
$\centra_G(T)$, which yields the desired result.
\end{proof}

\subsection{Residually finite lattices}

\begin{thm}\label{thm:ResiduallyFiniteCAT0Lattice}
Let $X$ be a proper \cat space such that $\Isom(X)$ acts cocompactly and minimally. Let $\Gamma < \Isom(X)$ be a
finitely generated  lattice. Assume that $\Gamma$ is irreducible
and residually finite. Then we have the following (see Section~\ref{sec:hull} for the notation):
\begin{enumerate}
\item  $\Gamma^*$ acts faithfully on each irreducible factor of $X$.

\item If $\Gamma$ is cocompact and $X$ is reducible, then for any closed subgroup $G < \Isom(X)$ containing
$H_{\Gamma^*}$, we have $\QZ(G)= \QZ(G^*)=1$. Furthermore $\soc(G^*)$ is a direct product of $r$ non-discrete
closed subgroups, each of which is characteristically simple, where $r$ is the number of irreducible factors of
$X$.
\end{enumerate}
\end{thm}

\noindent
(We emphasise that the irreducibility assumption concerns $\Gamma$ as an abstract group; compare however Remark~\ref{rem:hull}.)

\begin{proof}
If $X$ is irreducible, there is nothing to prove. We assume henceforth that $X$ is reducible. In view of
Theorem~\ref{thm:lattices:EuclideanSplitting}, $X$ has no Euclidean factor. Moreover,
Corollary~\ref{cor:Lattice=>NoFixedPointAtInfty} implies that $\Isom(X)$ fixes no point at infinity.
In particular, $\Gamma$ and $H_{\Gamma^*}$ act minimally without fixed point at infinity by
Theorem~\ref{thm:density}.

Let $H_1, \dots, H_r$ be the irreducible factors of $H_{\Gamma^*}$; thus $r$ coincides with the number of
irreducible factors of $X$. In view of Theorem~\ref{thm:lattices:Irred}, the group $\Gamma^*$ is an
irreducible lattice in this product.  By Corollary~\sref{cor:Isom(irreducible)} and Theorem~\ref{thm:density}, each $H_i$
is either a centre-free simple Lie group or totally disconnected with trivial amenable radical. If $H_1$ is a
simple Lie group, then it has no non-trivial discrete normal subgroup and hence
$$(\Gamma^*)_1 := \Gamma^* \cap (H_1 \times 1 \times \cdots \times 1) = 1.$$
If $H_1$ is totally disconnected, then by Proposition~\ref{prop:ResiduallyFinite} the canonical discrete
kernel $(\Gamma^*)_1 $ commutes with the discrete residual $H_1^{(\infty)}$, which is non-trivial by
Corollary~\sref{cor:ResiduallyDiscrete}. Thus $\centra_{H_1}(H_1^{(\infty)}) = 1$ by
Theorem~\sref{thm:geometric_simplicity} and hence $(\Gamma^*)_1  = 1$.

Assertion~(i) now follows from a straightforward induction on~$r$.

\medskip Assume next that $\Gamma$ is cocompact. Let $G_1, \dots, G_r$ be the irreducible factors of $G^*$. By
Lemma~\ref{lem:hull} and Proposition~\ref{prop:DiscreteNormal}, and in view of Part~(i), for each totally
disconnected factor $G_i$, the quasi-centre $\QZ(G_i)$ is topologically locally finite. Its closure is thus
amenable, hence trivial by Theorem~\sref{thm:geometric_simplicity}. Moreover, the quasi-centre of each almost
connected factor is trivial as well by Lemma~\ref{lem:hull}.

Clearly the projection of the quasi-centre of $G^*$ to the irreducible factor $G_i$ is contained in $\QZ(G_i)$.
This shows that $\QZ(G^*)$ is trivial. Hence so is $\QZ(G)$, since it contains $\QZ(G^*)$ as a finite index
subgroup and since $G$ has no non-trivial finite normal subgroup by Corollary~\sref{cor:NoEuclideanFactor}. Now
the desired conclusion follows from Proposition~\sref{prop:socle}.
\end{proof}

\subsection{Strong rigidity for product spaces}\label{sec:Mostow}
In~\cite{Caprace-Monod_structure}, we presented a few superrigidity results (Section~\sref{sec:few:super}).
Superrigidity should contain, in particular, strong rigidity \emph{à la} Mostow. This is indeed the content of
Theorem~\ref{thm:Mostow} below, where an isomorphism of lattices is shown to extend to an ambient group.
However, in contrast to the classical case of symmetric spaces, which are homogeneous, the full isometry group
does not in general determine the space since \cat spaces are in general not homogeneous. Another difference is
that the hull of a lattice, as described in Section~\ref{sec:hull}, is generally smaller than the full isometry
group of the ambient \cat space.

In view of the definition of the hull, the following statement is non-trivial only when $X$ (or an invariant subspace)
is reducible; this is expected since we want to use superrigidity for irreducible lattices in products.

\begin{thm}\label{thm:Mostow}
Let $X,Y$ be proper \cat spaces and $\Gamma, \Lambda$ discrete cocompact groups of isometries of $X$, respectively $Y$,
not splitting (virtually) a $\ZZ^n$ factor.

Then any isomorphism $\Gamma \cong \Lambda$ determines an isomorphism
$H_\Gamma\cong H_\Lambda$ such that the following commutes:
$$\xymatrix{\Gamma\ar[d]\ar[rr] && \Lambda\ar[d]\\
H_\Gamma\ar[rr]^{\cong} && H_\Lambda}$$
\end{thm}

Theorem~\ref{thm:Mostow} provides a partial answer to Question~21 in~\cite{FarbHruskaThomas}.

\begin{remark}
The assumption on $\ZZ^n$ factors is equivalent to excluding Euclidean factors from $X$ (or its canonical
invariant subspace) by Theorem~\ref{thm:lattices:EuclideanSplitting}. On the one hand, this assumption is really
necessary for the theorem to hold, even for symmetric spaces, since one can twist the product using a $\Gamma$-action
on the Euclidean factor when $\mathrm{H}^1(\Gamma)\neq 0$ (compare~\cite[\S\,4]{Lawson-Yau}).
On the other hand, since Bieberbach groups are obviously Mostow-rigid, Theorem~\ref{thm:Mostow} together with
Theorem~\ref{thm:lattices:EuclideanSplitting} give us as complete as possible a description of the situation
with $\ZZ^n$ factors.
\end{remark}

\begin{proof}
Let $X'\se X$ be the subset provided by Theorem~\ref{thm:Lattice=>NoFixedPointAtInfty}. We retain the notation
$F_\Gamma\lhd \Gamma$ and $\Gamma'=\Gamma/F_\Gamma < \Isom(X')$ from Section~\ref{sec:hull} and recall from
Lemma~\ref{lem:finiteRad} that $F_\Gamma$ depends only on $\Gamma$ as abstract group and that $X'$ is
$\Gamma'$-minimal. We define $Y'$, $F_\Lambda$ and $\Lambda'$ in the same way and have the corresponding lemma.
In particular, it follow that the isomorphism $\Gamma\cong\Lambda$ descends to $\Gamma'\cong\Lambda'$.
Therefore, we can and shall assume from now on that $X$ and $Y$ are minimal and $\Gamma<H_\Gamma<\Isom(X)$,
$\Lambda<H_\Lambda<\Isom(Y)$. By Theorem~\ref{thm:lattices:EuclideanSplitting}, we know that $X,Y$ have no
Euclidean factor. Thus $\Gamma, \Lambda$ have no fixed point at infinity by
Theorem~\ref{thm:Lattice=>NoFixedPointAtInfty_bis}. We claim that $\Gamma$ has a finite index subgroup
$\Gamma^\dagger$ which decomposes as a product $\Gamma^\dagger = \Gamma_1\times\cdots \times \Gamma_s$ of
irreducible factors, with $s$ maximal for this property. Indeed, otherwise we could apply the splitting theorem
of~\cite{Monod_superrigid} to a chain a finite index subgroups and contradict the properness of $X$. We write
$\Lambda^\dagger= \Lambda_1\times\cdots \times \Lambda_s$ for the corresponding groups in $\Lambda$. Combining
the splitting theorem with Addendum~\sref{addendum}, it follows from the definition of the hull that it is
sufficient to prove the statement for $s=1$. We assume henceforth that $\Gamma$, and hence also $\Lambda$, is
irreducible. Furthermore, if $X$ and $Y$ are both irreducible, then $H_\Gamma= \Gamma$ and $H_\Lambda = \Lambda$
and the desired statement is empty. We now assume that $X$ is reducible.

By Theorem~\ref{thm:density}, the lattice $\Gamma$ (resp. $\Lambda$) acts minimally without fixed point at
infinity on $X$ (resp. $Y$). Theorem~\sref{thm:Monod:superrigidity} yields a continuous morphism $f:H_{\Gamma^*} \to H_{\Lambda^*}$,
which shows in particular (by the splitting theorem~\cite{Monod_superrigid}) that $Y$ is reducible as well.
A second application of Theorem~\sref{thm:Monod:superrigidity} yields a second continuous  morphism $f':H_{\Lambda^*} \to
H_{\Gamma^*}$. Notice that the respective restrictions to $\Gamma^*$ and $\Lambda^*$ coincides with the given
isomorphism and its inverse. In particular $f' \circ f$ (resp. $f'\circ f$) is the identity on $\Gamma$ (resp.
$\Lambda$). By definition of the hull, it follows that $f' \circ f$ (resp. $f'\circ f$) is the
identity on $H_{\Gamma^*}$ (resp. $H_{\Lambda^*}$). The desired result finally follows, since there is a
canonical isomorphism $\Gamma/\Gamma^* \cong H_\Gamma/ H_{\Gamma^*}$ and since the action of $H_\Gamma/
H_{\Gamma^*}$ on $H_{\Gamma^*}$ is canonically determined by the action of $\Gamma/\Gamma^* $ on $\Gamma^*$.
\end{proof}

The above proof shows in particular that amongst spaces that are $\Gamma$-minimal without Euclidean factor, the
number of irreducible factors depends only upon the group $\Gamma$. If we combine this with
Theorem~\ref{thm:Lattice=>NoFixedPointAtInfty_bis}, Corollary~\ref{cor:lattices:EuclideanSplitting} and
Theorem~\ref{thm:lattices:EuclideanSplitting}\eqref{pt:lattices:EuclideanSplitting:min}, we obtain that the
number of factors in the ``de Rham'' decomposition
\begin{equation}\label{eq:deRham}
X' \cong\ X_1\times \cdots \times X_p \times \RR^n \times Y_1\times \cdots \times Y_q
\end{equation}
of Addendum~\sref{addendum} is an invariant of the group:

\begin{cor}\label{cor:deRham:lattice}
Let $X$ be a proper \cat space and $\Gamma < \Isom(X)$ be a group acting properly discontinuously and
cocompactly.

Then any other such space admitting a proper cocompact $\Gamma$-action has the same number of factors
in~\eqref{eq:deRham} and the Euclidean factor has same dimension.\qed
\end{cor}

\noindent
(We recall that minimality is automatic when $X$ is geodesically complete: Lemma~\sref{lem:cocompact:minimal}.)

\section{Arithmeticity of abstract lattices}\label{sec:arith}
The main goal of this section is to prove Theorem~\ref{thm:arith:general:intro}, which we now
state in a slightly more general form. Following G.~Margulis~\cite[IX.1.8]{Margulis}, we shall say that a
simple algebraic group $\mathbf{G}$ defined over a field $k$ is \textbf{admissible}\index{admissible group}
if \emph{none} of the following holds:

\begin{itemize}
\item[---] $\chr(k) = 2$ and $\mathbf{G}$ is of type $A_1$, $B_n$, $C_n$ or $F_4$.

\item[---] $\chr(k) = 3$ and $\mathbf{G}$ is of type $G_2$.
\end{itemize}

\noindent
A semi-simple group will be said admissible if all its factors are.

\begin{thm}\label{thm:arith:general}
Let $\Gamma<G=G_1\times\cdots\times G_n$ be an irreducible finitely generated lattice, where each $G_i$ is
any locally compact group.

If $\Gamma$ admits a faithful Zariski-dense representation in an admissible semi-simple group (over any field),
then the amenable radical $R$ of $G$ is compact and the quasi-centre $\QZ(G)$ is
virtually contained in $\Gamma\cdot R$. Furthermore, upon replacing $G$ by a finite index subgroup, the quotient
$G/R$ splits as $G^+ \times \QZ(G/R)$ where $G^+$ is a semi-simple algebraic group and the image of $\Gamma$ in
$G^+$ is an arithmetic lattice.
\end{thm}

Since the projection map $G \to G/R$ is proper, the statement of Theorem~\ref{thm:arith:general} implies in
particular that $\QZ(G/R)$ is discrete.

\begin{cor}\label{cor:arith:general}
Let $G=G_1\times\cdots\times G_n$ be a product of locally compact groups. Assume that $G$ admits a finitely
generated irreducible lattice with a faithful Zariski-dense representation in a semi-simple group over some
field of characteristic~$\neq 2,3$.

Then $G$ is a compact extension of a direct product of a semi-simple algebraic group by a discrete group.\qed
\end{cor}

To be more precise, the arithmeticity conclusion\index{arithmetic group}
of Theorem~\ref{thm:arith:general} means the following.

\smallskip

There exists a global field $K$, a connected semi-simple $K$-anisotropic $K$-group $\HH$ and a finite set
$\Sigma$ of valuations of $K$ such that:\label{page:arithmetic_group}

(i)~The quotient $\ud\Gamma:=\Gamma / \Gamma \cap (R \cdot \QZ(G))$ is commensurable with the arithmetic group
$\HH(K(\Sigma))$, where $K(\Sigma)$ is the ring of $\Sigma$-integers of $K$. Moreover, $\Sigma$ contains all
Archimedean valuations $v$ for which $\HH$ is $K_v$-isotropic, where $K_v$ denotes the $v$-completion of $K$. In
particular, by Borel--Harish-Chandra and Behr--Harder reduction theory, the diagonal embedding realises $\HH(K(\Sigma))$ as a
lattice in the product $\prod_{v\in \Sigma} \HH(K_v)$.

(ii)~The group $G^+$ is isomorphic to $\prod_{v\in \Sigma} \HH(K_v)^+$ and this isomorphism implements the
commensurability of $\ud\Gamma$ with $\HH(K(\Sigma))$. For background references, including on $\HH(K_v)^+$,
see~\cite[I.3]{Margulis}.

\medskip

In contrast to statements in~\cite{Monod_alternative}, there is no assumption on the subgroup structure of the
factors $G_i$ in Theorem~\ref{thm:arith:general}, which may not even be irreducible. The nature of the linear
representation is however more restricted.

Another improvement is that no (weak) cocompactness assumption is made on $\Gamma$. In particular, under the
same algebraic restrictions on the factors $G_i$ as in \emph{loc.\ cit.},
we obtain the following arithmeticity vs. non-linearity alternative for all finitely generated lattices.

\begin{cor}\label{cor:arith:SimpleFactors}
Let $\Gamma<G=G_1\times\cdots\times G_n$ be an irreducible finitely generated lattice, where each $G_i$ is a
locally compact group such that every non-trivial closed normal subgroup is cocompact. Then one of the following holds:

\begin{enumerate}
\item Every finite-dimensional linear representation of $\Gamma$ in characteristic~$\neq 2,3$ has virtually
soluble image.\label{pt:alternative:abelian}

\item $G$ is a semi-simple algebraic group and $\Gamma$ is an arithmetic lattice.\label{pt:alternative:arith}
\end{enumerate}
\end{cor}

The hypothesis made on each factor $G_i$ may be used to describe to some extent its structure independently of the existence of
a lattice in $G$; one can in particular show~\cite{Caprace-Monod_monolith} that each $G_i$ is monolithic, thus extending the
classical result of Wilson~\cite{Wilson71} to locally compact groups.
However, we will not appeal to this preliminary description of the $G_i$ when proving
Corollary~\ref{cor:arith:SimpleFactors}: the structural information will instead be obtained \emph{a posteriori}.

\begin{remark}\label{rem:superrigidity:char}
In~\cite{Monod_alternative}, the conclusion~\eqref{pt:alternative:abelian} was replaced by finiteness of the image.
This follows from the current conclusion in the more restricted setting of \emph{loc.\ cit.}
thanks to Y.~Shalom's superrigidity for characters~\cite{Shalom00}, unless of course $G_i$ admits (virtually) a non-zero
continuous homomorphism to $\RR$ (after all in the current setting we can have $G_i=\RR$). It is part of the
assumptions in~\cite{Monod_alternative} that no such homomorphism exists, so that Corollary~\ref{cor:arith:SimpleFactors}
indeed generalises  \emph{loc.\ cit}.
\end{remark}

\subsection{Superrigid pairs}
For convenience, we shall use the following terminology. Let $J$ be a topological group and $\Lambda < J$ any
subgroup. We call the pair $(\Lambda, J)$ \textbf{superrigid} if for any local field $k$ and any connected
absolutely almost simple adjoint $k$-group $\HH$, every abstract homomorphism $\Lambda \to \HH(k)$ with
unbounded Zariski-dense image extends to a continuous homomorphism of $J$.

\begin{prop}\label{prop:SuperrigidPairs}
Let $(\Lambda, J)$ be a superrigid pair with $J$ locally compact and $\Lambda$ finitely generated with closure
of finite covolume in $J$.

If $\Lambda$ admits a faithful representation in an admissible semi-simple group (over any field)
with Zariski-dense image, then the amenable radical $R$ of $J$ is compact and the quasi-centre $\QZ(J)$ is
virtually contained in $\Lambda\cdot R$. Furthermore, upon replacing $J$ by a finite index subgroup, the
quotient $J/R$ splits as $J^+ \times \QZ(J/R)$ where $J^+$ is a semi-simple algebraic group and the image of
$\Lambda$ in $J^+$ is an arithmetic lattice.
\end{prop}
\noindent
(We point out again that in particular the direct factor $\QZ(J/R)$ is discrete.)

\medskip

One might expect that Theorem~\ref{thm:arith:general} could now be proved by establishing in complete generality
that finitely generated irreducible lattices in products of locally compact groups form a superrigid pair. For
uniform lattices, or more generally weakly cocompact square-summable lattices, this is indeed true and was
proved in~\cite{Monod_superrigid}. We do not have a proof in general and shall eschew this difficulty by giving
first an independent proof of the compactness of the amenable radical (Corollary~\ref{cor:radical:rigid} below)
and using the residual finiteness of finitely generated linear groups before proceeding with
Proposition~\ref{prop:SuperrigidPairs}.

\smallskip

Nevertheless, we do have a general proof as soon as the groups are totally disconnected.

\begin{thm}\label{thm:SuperrigidPair:td}
Let $\Gamma<G=G_1\times\cdots\times G_n$ be an irreducible finitely generated lattice, where each $G_i$ is
any locally compact group.

If $G$ is totally disconnected, then $(\Gamma, G)$ is a superrigid pair.
\end{thm}

\noindent
(As we shall see in Proposition~\ref{prop:commensurators}, one can drop the finite
generation assumption in the simpler case where $\Gamma$ projects faithfully to some factor $G_i$.)

\smallskip

The proofs will use the following fact established in~\cite{Caprace-Monod_monolith}.

\begin{prop}\label{prop:SimpleQuotients}
Let $G$ be a compactly generated locally compact group and $\{N_v \; | \; v \in \Sigma\}$ be a collection of
pairwise distinct closed normal subgroups of $G$ such that for each $v \in \Sigma$, the quotient $H_v = G/N_v$
is quasi-simple, non-discrete and non-compact. If $\bigcap_{v \in \Sigma} N_v = 1$ then $\Sigma$ is finite and
$G$ has a characteristic closed cocompact subgroup which splits as a finite direct product of $|\Sigma|$
topologically simple groups.\qed
\end{prop}

We recall for the above statement that a group is called \textbf{quasi-simple}\label{def:quasi-simple}\index{quasi-simple} if it
possesses a cocompact normal subgroup which is topologically simple and contained in every non-trivial closed
normal subgroup.

\begin{proof}[Proof of Proposition~\ref{prop:SuperrigidPairs}]
We will largely follow the ideas of Margulis, deducing arithmeticity from
superrigidity~\cite[Chapter~IX]{Margulis}. It is assumed that the reader has a copy of~\cite{Monod_alternative}
at hand, since it contains a similar reasoning under different hypotheses. The characteristic assumption in \emph{loc.\ cit.}
will be replaced by the current admissibility assumption.

The group $J$ (and hence also all finite index subgroups and factors) is compactly generated
by Lemma~\ref{lem:UGammaU}. Let $\tau : \Lambda \to \HH$ be the given faithful representation.
 Upon replacing $\Lambda$ and $J$ by finite index subgroups and post-composing $\tau$ with the projection map $\HH \to
\HH/\centra(\HH)$, we shall assume henceforth that $\HH$ is adjoint and Zariski-connected. The representation
$\tau: \Lambda \to \HH$ need no longer be faithful, but it still has finite kernel.
As in~\cite[(3.3)]{Monod_alternative}, in view of the assumption that $\Lambda$ is finitely generated,
we may assume that $\HH$ is defined over a finitely generated field $K$. This is the first of two places where
the admissibility assumption is used in \emph{loc.\ cit.} (following~VIII.3.22 and~IX.1.8 in~\cite{Margulis}).

By Tits' alternative~\cite{Tits72}, the amenable radical of $\Lambda$ is soluble-by-locally-finite
and thus locally finite since $\tau(\Lambda)$ is Zariski-dense and $\HH$ is semi-simple.
The finite generation of $K$ implies that this radical is in fact finite (see \emph{e.g.} Corollary~4.8 in~\cite{Wehrfritz}),
thus trivial by Zariski-density since
$\HH$ is adjoint. (The finite generation of $K$ is essential in positive characteristic since for algebraically
closed fields there is always a locally finite \emph{Zariski-dense} subgroup~\cite{Bernik-Guralnick-Mastnak}.)

It now follows that if $J$ is a compact extension of a discrete group, then the latter has trivial amenable radical
and thus all the conclusions of Proposition~\ref{prop:SuperrigidPairs} hold trivially. Therefore,
we assume henceforward that $J$ is not compact-by-discrete.

\smallskip
Let $\HH = \HH_1 \times \dots \times \HH_k$ be the decomposition of $\HH$ into its simple factors.
We shall work with the factors $\HH_i$ one at a time. Let $\tau_i : \Lambda \to \HH_i$ be the induced representation of
$\Lambda$. Notice that $\tau_i$ need not be faithful; however, it has Zariski-dense (and in particular infinite)
image.

We let $\Sigma_i$ denote the set of all (inequivalent representatives of) valuations $v$ of $K$ such that the image of
$\tau_i(\Lambda)$ is not relatively compact in $\HH_i(K_v)$ (for the Hausdorff topology); observe that this
image is still Zariski-dense. Then $\Sigma_i$ is non-empty since $\tau_i(\Lambda)$ is infinite, see~\cite[Lemma~2.1]{Breuillard-Gelander}.

By hypothesis, there exists a continuous representation $J \to \HH_i(K_v)$ for each $v \in \Sigma_i$, extending
the given $\Lambda$-representation. We denote by $N_v\lhd J$ the kernel of this representation. Let $I \subseteq
\{1, \dots, k\}$ be the set of all those indices $i$ such that $J/N_v$ is non-discrete for each $v \in
\Sigma_i$.

\smallskip
\noindent \emph{We claim that the set $I$ is non-empty.}

\smallskip

Indeed, for each index $j \not \in I$, there exists $v_j \in \Sigma_j$ such that $N_{v_j}$ is open in $J$. Thus
the kernel
$$J^+ = \bigcap_{j \not \in I} N_{v_j} $$
of the continuous representation $J \to \prod_{j \not \in I}\HH_j(K_{v_j})$ is open.

By assumption the closure of $\Lambda$ in $J$ has finite covolume. Therefore, for each open subgroup $F < J$,
the closure of $\Lambda \cap F$ has finite covolume in $F$.
It follows in particular that $\Lambda \cap F$ is infinite unless $F$ is compact.

These considerations apply to the open subgroup $J^+ < J$. Since $J$ is not compact-by-discrete, we deduce that
$\Lambda \cap J^+$ is infinite. Therefore the restriction to $\Lambda$ of the representation $J \to \prod_{j
\not \in I}\HH_j(K_{v_j})$ has infinite kernel and, hence, it does not factor through $\tau : \Lambda \to
\HH(K)$. In particular it cannot coincide with the given representation $\tau : \Lambda \to \HH$. Thus $I$ is
non-empty.

\smallskip
\noindent \emph{We claim that for each $i \in I$, the set $\Sigma_i$ is finite.}

\smallskip

Let $i \in I$ and $v \in \Sigma_i$. The arguments of~\cite[(3.7)]{Monod_alternative} show that the isomorphic
image of $J/N_v$ in $\HH_i(K_v)$ contains $\HH_i(K_v)^+$. These arguments use again the admissibility assumption
because the appeal to a result of R.~Pink~\cite{Pink98}; the fact that the latter hold in the admissible case
is explicit in the table provided in Proposition~1.6 of~\cite{Pink98}. Furthermore, it follows from Tits' simplicity theorem~\cite{Ti64}
combined with~\cite[6.14]{Borel-Tits73} that each $J/N_v$ is quasi-simple.
Moreover, an application of~\cite[8.13]{Borel-Tits73}  shows
that the various quotients $(J/N_v)_{v \in \Sigma_i}$ are pairwise non-isomorphic. In particular the normal
subgroups $(N_v )_{v \in \Sigma_i}$ are pairwise distinct.

Let $D_i = \bigcap_{v \in \Sigma_i} N_v$ and recall that $J/D_i$ is compactly generated. Projecting each $N_v$
to $J/D_i$, we obtain a family of pairwise distinct normal subgroups of $J/D_i$ indexed by $\Sigma_i$ such that
each corresponding quotient is quasi-simple, non-discrete and non-compact. Therefore, the desired claim follows
from Proposition~\ref{prop:SimpleQuotients}.

\smallskip

In particular, appealing again to~\cite[Corollaire~8.13]{Borel-Tits73},  we obtain a continuous map $J \to
\prod_{v \in \Sigma_i} \HH_i(K_v)$ which we denote again by $\tau^J_i$. The kernel of $\tau^J_i$ is $D_i$.  Upon
replacing $J$ and $\Lambda$ by finite index subgroups we may assume that the image of $\tau^J_i$ coincides in
fact with $\prod_{v\in \Sigma_i} \HH_i(K_v)^+$, compare~\cite[(3.9)]{Monod_alternative}.

\smallskip
\noindent \emph{We claim that $R := J^+ \cap D$ is compact and that $J =  J^+ \cdot D$, where $D$ is defined by
$D = \bigcap_{i \in I} D_i$.}

\smallskip

We first show that $R =J^+ \cap D$ is compact. Assume for a contradiction that this is not the case.  Then,
given a compact open subgroup $U$ of $J$, the intersection $\Lambda_0 = \Lambda \cap (U \cdot R)$ is infinite:
this follows from the same argument as above, using the assumption that the closure of $\Lambda$ has finite
covolume.

For each index $j \not \in I$, we have  $J^+ \se N_{v_j}$ and we deduce that the image of $\Lambda_0$ in
$\HH_j(K_{v_j})$ is finite, since it is contained in the image of $U$. Equivalently, the subgroup
$\tau_j(\Lambda_0) < \HH_j(K)$ is finite. It follows in particular that $\tau_i(\Lambda_0)$ is infinite for some
$i \in I$. By~\cite[Lemma~2.1]{Breuillard-Gelander}, there exists $v \in \Sigma_i$ such that the image of
$\Lambda_0$ in $\HH_i(K_v)$ is unbounded. This is absurd since $D \se N_v$ and hence the image of $\Lambda_0$ in
$\HH_i(K_v)$ is contained in the image of the compact subgroup $U$. This shows that the intersection $R$ is
indeed compact.

At this point we know that the quotient $J/D$ is isomorphic to a subgroup of the product
$$\prod_{i \in I}\prod_{ v\in \Sigma_i} \HH_i(K_v)^+$$
which projects surjectively onto each factor of the form $\prod_{ v\in \Sigma_i} \HH_i(K_v)^+$.
Using again the Goursat-type argument as in Proposition~\ref{prop:SimpleQuotients},
we find that $J/D$ is indeed isomorphic to a finite product of non-compact non-discrete simple groups
$\HH_i(K_v)^+$. In particular the quotient $J/D$ has no non-trivial open normal subgroup. Since $J^+$ is open and normal in $J$,
we deduce that $J = J^+\cdot D$, thereby establishing the claim.

\medskip

By the very nature of the statement, we may replace $J$ by the quotient $J/R$ without any loss of generality,
since $R$ is compact. In view of this further simplification, the preceding claim implies that $J \cong J^+
\times D$. In particular $D$ is discrete.

It now follows as in~\cite[(3.11)]{Monod_alternative} that $K$ is a global field, and that the image of $\Lambda$
in the semi-simple group $J/D$ is an arithmetic lattice (compare~\cite[(3.13)]{Monod_alternative}). Therefore, by
Proposition~\ref{prop:Raghunathan}, the intersection $\Lambda \cap D$ is a lattice in $D$ and, hence, the
discrete normal subgroup $D$ is virtually contained in $\Lambda$. As $J \cong J^+ \times D$ and $J^+$ has trivial
quasi-centre, it follows that the quasi-centre of $J$ coincides with $D$. This finishes the proof.
\end{proof}

For later use, we single out a (simpler) version of an argument referred to above.

\begin{lem}\label{lem:containsH+}
Let $\HH$ be an admissible connected absolutely almost simple adjoint $k$-group $\HH$, where $k$ is a local field.
Let $J$ be a locally compact group with a continuous unbounded Zariski-dense homomorphism $\tau : J\to\HH(k)$.
Then any compact normal subgroup of $J$ is contained in the kernel of $\tau$.
\end{lem}

\begin{proof}
Let $K\lhd J$ be a compact normal subgroup. The Zariski closure of $\tau(K)$ is normalised by the Zariski-dense group
$\tau(J)$ and therefore it is either $\HH(k)$ or trivial. We assume the former since otherwise we are done.

We claim that we can assume $k$ non-Archimedean. Otherwise, either $k=\RR$ or $k=\CC$. In the first case, $\tau(K)$ coincides with its Zariski
closure by Weyl's algebraicity theorem~\cite[4.2.1]{VinbergENC} so that $\HH(k)$ is compact in which case the lemma is void by the unboundedness
assumption. In the second case, one can reduce to the case $\tau(J)\se\HH(\RR)$ as in~\cite[(3.5)]{Monod_alternative} and thus $\tau(K)=1$ as before
since $\HH(\RR)$ is also simple; the claim is proved.

Following now an idea from~\cite[p.~41]{Shalom00} (see also the explanations in Section~(3.7)
of~\cite{Monod_alternative}), one uses~\cite{Pink98} to deduce that $\tau(K)$ is open upon possibly replacing
$k$ by a closed subfield (the admissibility assumption enters as in the proof of Proposition~\ref{prop:SuperrigidPairs}).
We can still denote this subfield by $k$ because it accommodates the whole image
$\tau(J)$, see again~\cite[(3.7)]{Monod_alternative}. Now $\tau(J)$ is an unbounded open subgroup and hence
contains $\HH(k)^+$ by a result of J.~Tits (see~\cite{Prasad82}; this also follows from the Howe--Moore
theorem~\cite{Howe-Moore} which however is posterior to Tits' result). This implies that the compact group
$\tau(K)$ is trivial since $\HH(k)^+$ is simple by~\cite{Ti64}.
\end{proof}

\subsection{Boundary maps}
We record two statements extracted from Margulis' work in the form most convenient for us.

\begin{prop}\label{prop:furst:extension}
Let $J$ be a second countable locally compact group with a measure class preserving action on a standard
probability space $B$. Let $\Lambda < J$ be a dense subgroup with a Zariski-dense unbounded representation
$\tau:\Gamma\to \HH(k)$ to a connected absolutely almost $k$-simple adjoint group $\HH$ over an arbitrary local field $k$.

If there is a proper $k$-subgroup $\LL<\HH$ and a $\Lambda$-equivariant non-essentially-constant measurable map
$B\to \HH(k)/\LL(k)$, then $\tau$ extends to a continuous homomorphism $J \to \HH(k)$.
\end{prop}

\begin{proof}
The argument is given by A'Campo--Burger in the characteristic zero case at the end of Section~7
in~\cite{A'Campo-Burger} (pp.~18--19). This reference considers homogeneous spaces for $B$ but this restriction
is never used. The general statement is referred to in~\cite{Burger_ICM} and details are given in~\cite{Bonvin}.
\end{proof}

\begin{prop}\label{prop:furst:exists}
Let $\Gamma$ be a countable group with a Zariski-dense unbounded representation $\Gamma\to \HH(k)$ to a
connected absolutely almost $k$-simple adjoint group $\HH$ over an arbitrary local field $k$.
Let $B$ be a standard probability space with a measure class preserving $\Gamma$-action that is amenable in Zimmer's sense~\cite{Zimmer84}
and such that the diagonal action on $B^2$ is ergodic.

Then there is a proper $k$-subgroup $\LL<\HH$ and a $\Gamma$-equivariant non-essentially-constant measurable map
$B\to \HH(k)/\LL(k)$.
\end{prop}

\begin{proof}
Again, this is proved in~\cite{A'Campo-Burger} for the characteristic zero case (and $B$ homogeneous) and the necessary adaptations
to the general case are explained in~\cite{Bonvin}.
\end{proof}

We shall need these specific statements below. They first appeared within the proof of Margulis' commensurator superrigidity, which can
adapted as follows using~\cite{Burger_ICM} and Lemma~\sref{lem:elementary:superrigidity}, providing a first step towards
Theorem~\ref{thm:SuperrigidPair:td}.

\begin{prop}\label{prop:commensurators}
Let $G=G_1\times G_2$ be a product of locally compact $\sigma$-compact groups and $\Lambda<G$ be an irreducible lattice.
Assume that the projection of $\Lambda$ to $G_1$ is injective and that
$G_2$ admits a compact open subgroup. Then the pair $(\Lambda, G)$ is superrigid.
\end{prop}

\begin{proof}
We claim that one can assume $G$ second countable.
As explained in~\cite[Proposition~61]{Monod_superrigid}, $\sigma$-compactness implies the existence
of a compact normal subgroup $K\lhd G$ meeting $\Lambda$ trivially and such that $G/K$ is second countable. Applying the
statement to $G/K$ together with the image of $\Lambda$ therein yields the general statement since the projection of $\Lambda$
to $G/K$ is an isomorphism; this proves the claim.

\smallskip
Let $\tau : \Lambda \to \HH(k)$ be as in the definition of superrigid pairs
and let $U < G_2$ be a compact open subgroup. Set $\Lambda_U  = \Lambda \cap (G_1 \times U)$.
By the injectivity assumption and Lemma~\ref{lem:CompactOpenTrick}, we can consider $\Lambda_U$
as a lattice in $G_1$ which is commensurated by (the image of the projection of) $\Lambda$. We distinguish two cases.

Assume first that $\tau(\Lambda_U)$ is unbounded in the locally compact group $\HH(k)$. We may then apply
Margulis' commensurator superrigidity in its general form proposed by M.~Burger~\cite[Theorem~2.A]{Burger_ICM},
see~\cite{Bonvin} for details. This yields a continuous map $J \to \HH(k)$ factoring through $G_1$ and extending
the given $\Lambda$-representation, as desired.

Assume now that $\tau(\Lambda_U)$ is bounded, which is equivalent to $\Lambda_U$ fixing a point in the symmetric
space or Bruhat--Tits building associated to $\HH(k)$. Then Lemma~\sref{lem:elementary:superrigidity}
yields a continuous map $J \to \HH(k)$ factoring through $G_2$.
\end{proof}

\subsection{Radical superrigidity}

\begin{thm}\label{thm:radical:rigid}\index{radical!superrigidity}
Let $G$ be a locally compact group, $R\lhd G$ its amenable
radical, $\Gamma<G$ a finitely generated lattice and $F$ the closure of the image of $\Gamma$ in $G/R$.

Then any Zariski-dense unbounded representation of $\Gamma$ in any connected absolutely almost simple adjoint
$k$-group $\HH$ over any local field $k$ arises from a continuous representation of $F$ \emph{via} the map $\Gamma\to F$.
\end{thm}

\noindent
(In particular, the pair $(\Gamma/(\Gamma\cap R), F)$ is superrigid.)

\begin{proof}
Notice that $G$ is $\sigma$-compact since it contains a finitely generated, hence countable, lattice. (In fact
$G$ is even compactly generated by Lemma~\ref{lem:UGammaU}.) Set $J=G/R$. There exists a standard probability
$J$-space $B$ on which the $\Gamma$-action is amenable and such that the diagonal $\Gamma$-action on $B^2$ is
ergodic; it suffices to choose $B$ to be the Poisson boundary of a symmetric random walk with full support on
$J$. Indeed: (i)~The $J$-action is amenable as was shown by Zimmer~\cite{Zimmer78b}; this implies that the
$G$-action is amenable since $R$ is an amenable group and thus that the $\Gamma$-action is amenable since
$\Gamma$ is closed in $G$ (see~\cite[5.3.5]{Zimmer84}). (ii)~The diagonal action of any closed finite covolume
subgroup $F<J$ on $B^2$ is ergodic in view of the \emph{ergodicity with coefficients} of $J$, and hence the same
holds for dense subgroups of $F$. For detailed background on this strengthening of ergodicity introduced
in~\cite{Burger-Monod} and on the Poisson boundary in general, we refer the reader to~\cite{Kaimanovich03}.

Let now $k$ be a local field, $\HH$ a connected absolutely almost simple $k$-group and $\Gamma\to \HH(k)$ a
Zariski-dense unbounded representation. We can apply Proposition~\ref{prop:furst:exists} and obtain a proper subgroup
$\LL<\HH$ and a $\Gamma$-equivariant map $B\to \HH(k)/\LL(k)$. Writing $\Lambda$ for the image of $\Gamma$ in $J$,
we can therefore apply Proposition~\ref{prop:furst:extension} with $F$ instead of $J$ and the conclusion follows.
\end{proof}

\begin{remark}
An examination of this proof shows that one has also the following related result. Let $J$ be a second countable
locally compact group and $\Lambda \se J$ a dense countable subgroup whose action on $J$ by left multiplication
is amenable. Then the pair $(\Lambda, J)$ is superrigid. Indeed, one can again argue with
Propositions~\ref{prop:furst:extension} and~\ref{prop:furst:exists} because it is easy to check that in the
present situation any amenable $J$-space is also amenable for $\Lambda$ viewed as a discrete group.
Related ideas were used by R.~Zimmer in~\cite{Zimmer87}.
\end{remark}

\begin{cor}\label{cor:radical:rigid}
Let $G$ be a locally compact group and $\Gamma<G$ a finitely generated lattice.

If $\Gamma$ admits a faithful Zariski-dense representation in an admissible semi-simple group (over any field),
then the amenable radical of $G$ is compact.
\end{cor}

\begin{proof}
Let $R$ be the amenable radical of $G$, $F$ be the closure of the image of $\Gamma$ in $G/R$ and $J < G$ the
preimage of $F$ in $G$. The content of Theorem~\ref{thm:radical:rigid} is that the pair $(\Gamma, J)$ is
superrigid. Since in addition $\Gamma$ is closed and of finite covolume in $J$ (see
\cite[Lemma~1.6]{Raghunathan}), we may apply Proposition~\ref{prop:SuperrigidPairs} and deduce that the amenable
radical of $J$ is compact. The conclusion follows since $R < J$.
\end{proof}

\subsection{Lattices with non-discrete commensurators}\label{sec:arith:commensurator}
The following useful trick allows to realise the commensurator of any lattice in a locally compact group $G$ as
a lattice in a product $G \times D$.
A similar reasoning in the special case of automorphism groups of trees may be found in~\cite[Theorem~6.6]{BenakliGlasner}.

\begin{lem}\label{lem:CommensuratorTrick}
Let $\Lambda$ be a group and $\Gamma<\Lambda$ a subgroup commensurated by $\Lambda$.
Let $D$ be the completion of $\Lambda$ with respect to the left or right uniform structure generated by the
$\Lambda$-conjugates of $\Gamma$. Then $D$ is a totally disconnected locally compact group.

If furthermore $G$ is a locally compact group containing $\Lambda$ as a dense subgroup such that
$\Gamma$ is discrete (resp. is a lattice) in $G$, then the diagonal embedding of $\Lambda$ in
$G \times D$ is discrete (resp. is an irreducible lattice).
\end{lem}

The above lemma is in some sense a converse to Lemma~\ref{lem:CompactOpenTrick}. In the special case where
one starts with a lattice satisfying a faithfulness condition, this relation becomes even stronger.

\begin{lem}\label{lem:BothTricks}
Let $G, H$ be locally compact groups and $\Lambda<G\times H$ a lattice. Assume that the projection of $\Lambda$ to
$G$ is faithful and that both projections are dense. Let $U<H$ be a compact open subgroup, set
$\Gamma=\Lambda \cap(G\times U)$ as in Lemma~\ref{lem:CompactOpenTrick} and consider the group $D$ as in
Lemma~\ref{lem:CommensuratorTrick} (upon viewing $\Lambda$ as a subgroup of $G$).
Define the compact normal subgroup $K\lhd H$ as the core $K=\bigcap_{h\in H} hUh\inv$ of $U$ in $G$.

Then the map $\Lambda\to D$ induces an isomorphism of topological groups $H/K \cong D$.
\end{lem}

\begin{proof}[Proof of Lemma~\ref{lem:CommensuratorTrick}]
One verifies readily the condition given in~\cite{BourbakiTGIII} (TG~III, \S\,3, No~4, Théorème~1) ensuring
that the completion satisfies the axioms of a group topology. We emphasise that it is part of
the definition of the completion that $D$ is Hausdorff; in other words $D$ is obtained by first completing
$\Lambda$ with respect to the group topology as defined above, and then dividing out the normal subgroup
consisting of those elements which are not separated from the identity.

Let $U$ denote the closure of the projection of $\Gamma$ to $D$. By definition $U$ is open. Notice that it is
compact since it is a quotient of the profinite completion of $\Gamma$ by construction. In particular $D$ is
locally compact.

By a slight abuse of notation, let us identify $\Gamma$ and $\Lambda$ with their images in $D$. We claim that $U
\cap \Lambda = \Gamma$. Indeed, let $\{\gamma_n\}_{n \geq 0}$ be a sequence of elements of $\Gamma$ such that
$\lim_n \gamma_n = \lambda \in \Lambda$. Since $\lambda \Gamma \lambda\inv$ is a neighbourhood of the identity
in $\Lambda$ (with respect to the topology induced from $D$), it follows that $\gamma_n \lambda\inv \in \lambda
\Gamma \lambda\inv$ for $n$ large enough. Thus $\lambda \in \gamma_n \Gamma = \Gamma$.

Assume now that $\Gamma$ is discrete and choose a neighbourhood $V$ of the identity in $G$ such that $\Gamma
\cap V = 1$. In view of the preceding claim the product $V \times U$ is a neighbourhood of the identity in $G
\times D$ which meets $\Lambda$ trivially, thereby showing that $\Lambda$ is discrete.

Assume finally that $\Gamma$ is a lattice in $G$ and let $F$ be a fundamental domain. Then $F \times U$ is a
fundamental domain for $\Lambda$ in $G \times D$, which has finite volume since a Haar measure for $G \times D$
may be obtained by taking the product of respective Haar measures for $G$ and $D$. Thus $\Lambda$ has finite
covolume in $G \times D$.
\end{proof}

\begin{proof}[Proof of Lemma~\ref{lem:BothTricks}]
In order to construct a continuous homomorphism $\pi:H\to D$, it suffices to check that any net in $\Lambda$ whose image in
$H$ converges to the identity also converges to the identity in $D$; this follows from the definitions of $\Gamma$ and $D$ since
the net is eventually in any $\Lambda$-conjugate of $U$. Notice that $\pi$ has dense image.

We claim that the kernel of $\pi$ is $\bigcap_{\lambda\in \Lambda} \lambda U\lambda\inv$.
Indeed, if on the one hand $k\in\ker(\pi)$ is the limit of the images in $H$ of a net $\{\lambda_i\}$ in $\Lambda$, then for any $\lambda$
we have eventually $\lambda_i\in \lambda\inv\Gamma\lambda \se \lambda\inv(G\times U)\lambda$ so that indeed $k\in \lambda\inv U\lambda$
since $U$ is closed. Conversely, if $k\in \bigcap_{\lambda\in \Lambda} \lambda U\lambda\inv$ is limit of images of $\{\lambda_i\}$,
then, since $U$ is open, for any $\lambda$ the image of $\lambda_i$ is eventually in $\lambda U\lambda\inv$, hence in
$\lambda \Gamma\lambda\inv$ so that $\pi(\lambda_i)\to 1$. This proves the claim.

Now it follows that $\ker(\pi)$ is indeed the core $K$ of the statement since $U$ is compact.
The fact that $\pi$ is onto and open follows from the existence of a compact open subgroup in $H$.
\end{proof}

\begin{thm}\label{thm:arith:commensurator}
Let $G$ be a locally compact group and $\Gamma<G$ be a lattice. Assume that $G$ possesses a
finitely generated dense subgroup $\Lambda$ such that $\Gamma < \Lambda < \Comm_G(\Gamma)$.

If $\Lambda$ admits a faithful Zariski-dense representation in an admissible semi-simple group (over any field),
then the amenable radical $R$ of $G$ is compact and the quasi-centre $\QZ(G)$ is
virtually contained in $\Gamma\cdot R$. Furthermore, upon replacing $G$ by a finite index subgroup, the quotient
$G/R$ splits as $G^+ \times \QZ(G/R)$ where $G^+$ is a semi-simple algebraic group and the image of $\Gamma$ in
$G^+$ is an arithmetic lattice.
\end{thm}

\begin{proof}
Let $J = G \times D$, where $D$ is the totally disconnected locally compact group provided by
Lemma~\ref{lem:CommensuratorTrick}. As a totally disconnected group, it has numerous compact open subgroups
(for instance the closure of $\Gamma$).
We shall view $\Lambda$ as an irreducible lattice in $J$. The projection of $\Lambda$ to $G$ is faithful by construction.
By Proposition~\ref{prop:commensurators}, the pair $(\Lambda, J)$ is superrigid.
This allows us to apply Proposition~\ref{prop:SuperrigidPairs}. Since the amenable
radical $R_G$ of $G$ is contained in the amenable radical $R_J$ of $J$, it is compact.
Furthermore, the quasi-centre of $G$ is contained in the quasi-centre of $J$ and the centre-free group $G/R_G$
is a direct factor of $J^+\times \QZ(J/R_J)$; the desired conclusions follow.
\end{proof}

\subsection{Lattices in products of Lie and totally disconnected groups}\label{sec:arith_statement}

\begin{thm}\label{thm:arith:SxD}
Let $\Gamma<G=S\times D$ be a finitely generated irreducible lattice, where $S$ is a connected semi-simple Lie
group with trivial centre and $D$ is a totally disconnected locally compact group. Let
$\Gamma_{\!D}\lhd D$ be the canonical discrete kernel of $D$.

Then $D/\Gamma_{\!D}$ is a profinite extension of a semi-simple algebraic group $Q$ and the image of $\Gamma$ in
$S\times Q$, which is isomorphic to $\Gamma/\Gamma_{\! D}$, is an arithmetic lattice.
\end{thm}

\begin{cor}\label{cor:arith:SxD}
In particular, $D$ is locally profinite by analytic.\hfill\qedsymbol
\end{cor}

A family of examples will be constructed in Section~\ref{sec:example} below, showing that the statement cannot be
simplified even in a geometric setting (see Remark~\ref{rem:profinite_extension}).

\begin{proof}[Proof of Theorem~\ref{thm:arith:SxD}]
By the very nature of the statement, we can factor out the canonical discrete kernel. Therefore, we shall assume
henceforth that the projection map $\Gamma\to S$ is injective. We can also assume that $S$ has no compact factors.
Since $S$ is connected with trivial centre, there is a Zariski connected semi-simple adjoint $\RR$-group $\HH$
without $\RR$-anisotropic factors such that $S=\HH(\RR)$. Notice that the injectivity of $\Gamma\to S$ is
preserved when passing to finite index subgroups.

By Proposition~\ref{prop:commensurators}, the pair $(\Gamma, G)$ is superrigid. We can therefore apply
Proposition~\ref{prop:SuperrigidPairs}. In particular, $D$ has compact amenable radical and therefore, in view
of the statement of Theorem~\ref{thm:arith:SxD}, we can assume that this radical is trivial. Given the
conclusion of Proposition~\ref{prop:SuperrigidPairs}, it only remains to show that the quasi-centre $\QZ(G)$ of
$G$ is trivial. We now know that $\QZ(G)$ is virtually contained in $\Gamma$; since on the other hand $S$ has
trivial quasi-centre, $\QZ(G)\se 1\times D$. In other words, $\QZ(G)$ is contained in the discrete kernel
$\Gamma_D$, which has been rendered trivial. This completes the proof.
\end{proof}

We have treated Theorem~\ref{thm:arith:SxD} as a port of call on the way to Theorem~\ref{thm:arith:general}.
In fact, one can also describe lattices in products of groups with a simple algebraic factor over an arbitrary local
field and in most cases without assuming finite generation \emph{a priori}.
We record the following statement, which will not be used below.

\begin{thm}\label{thm:MixedProducts}
Let $k$ be any local field and $\GG$ an admissible connected absolutely almost simple adjoint $k$-group.
Let $H$ be any compactly generated locally compact group admitting a compact open subgroup.
Let $\Gamma<\GG(k)\times H$ be an irreducible lattice.
In case $k$ has positive characteristic and the $k$-rank of $\GG$
is one, we assume $\Gamma$ cocompact.

Then $H/\Gamma_{\!H}$ is a compact extension of a semi-simple algebraic group $Q$ and the image of $\Gamma$ in
$\GG(k)\times Q$ is an arithmetic lattice.
\end{thm}

There is no assumption whatsoever on the compactly generated locally compact group $H$ beyond admitting a compact open subgroup;
recall that the latter is automatic if $H$ is totally disconnected~\cite[III \S\,4 No~6]{BourbakiTGI}. Notice that \emph{a posteriori} it follows from
arithmeticity that $\Gamma$ is finitely generated; in the proof below, finite generation will be established in two steps.


\begin{proof}[Proof of Theorem~\ref{thm:MixedProducts}]
We factor out the canonical discrete kernel $\Gamma_{\!H}$ and assume henceforth that it is trivial. This does
not affect the other assumptions and thus we choose some compact open subgroup $U<H$. We write $G=\GG(k)$ and
consider $\Gamma_U= \Gamma\cap(G\times U)$ as in Lemma~\ref{lem:CompactOpenTrick}. Since we factored out the
canonical discrete kernel, we can consider $\Gamma_U$ as a lattice in $G$ commensurated by the dense subgroup
$\Gamma<G$. Moreover, $\Gamma_U$ is finitely generated; indeed, either we have simultaneously $\rank_k(\GG)=1$
and $\chr(k)>0$, in which case we assumed $\Gamma$ cocompact, so that $\Gamma_U$ is cocompact in the compactly
generated group $\GG(k)$ (again Lemma~\ref{lem:CompactOpenTrick}) and hence finitely
generated~\cite[I.0.40]{Margulis}; or else, $\Gamma_U$ is known to be finitely generated by applying, as the
case may be, either Kazhdan's property, or the theory of fundamental domains, or the cocompactness of $p$-adic
lattices~--- we refer to Margulis, Sections~(3.1) and~(3.2)  of Chapter~IX in~\cite{Margulis}.

We can now apply Margulis' arthmeticity~\cite[1.(1)]{Margulis} and deduce that $\GG$ is defined over a global
field $K$ and that $\Gamma_U$ is commensurable to $\GG(K(S))$ for some finite set of places $S$; in short
$\Gamma_U$ is $S$-arithmetic. (The idea to obtain first this preliminary arithmeticity of $\Gamma_U$ was
suggested by M.~Burger.) It follows that $\Gamma$ is rational over the global field $K$, see
Theorem~3.b in~\cite{Borel66} (\emph{loc.\ cit.} is formulated for the Lie group case; see~\cite[Lemma~7.3]{Wortman}
in general).

\smallskip

Since the pair $(\Gamma, G\times H)$ is superrigid (for instance by Proposition~\ref{prop:commensurators}), only
the \emph{a priori} lack of finite generation for $\Gamma$ prevents us from applying
Proposition~\ref{prop:SuperrigidPairs}. However, a good part of the proof of that proposition is already secured
here since $\Gamma$ has been shown to be rational over a global field. We now proceed to explain how to adapt
the remaining part of that proof to the current setting. We use those elements of notation introduced
in the proof of Proposition~\ref{prop:SuperrigidPairs} that do not conflict with present notation and review all uses
of finite generation that are either explicit in the proof of Proposition~\ref{prop:SuperrigidPairs} or implicit
through references to~\cite{Monod_alternative}.

The compact generation of $G\times H$ is an assumption rather than a consequence of Lemma~\ref{lem:UGammaU}.

We also used finite generation in order to pass to a finite index subgroup of $\Gamma$ contained in $\GG(K_v)^+$ for
all valuations $v\in\Sigma$. We shall postpone this step, so that the whole argumentation provides us with maps
from $G\times H$ to a product $Q$ of factors that lie in-between $\GG(K_v)^+$ and $\GG(K_v)$. In particular all these
factors are quasi-simple and we can still appeal to Proposition~\ref{prop:SimpleQuotients} as before.
Notice however that at the very end of the proof, once finite generation is granted, we can invoke the argument
that $\GG(K_v)/\GG(K_v)^+$ is virtually torsion Abelian~\cite[6.14]{Borel-Tits73} and thus reduce again to the case
where $\Gamma$ is contained in $\GG(K_v)^+$.

We now justify that the image of $\Gamma$ in $G\times Q$ is discrete because previously this followed
from~\cite[(3.13)]{Monod_alternative} which relies on finite generation. If $\Gamma$ were not discrete, an
application of~\cite[Lemma~2.1]{Breuillard-Gelander} would provide a valuation $v\notin \Sigma$ with $\Gamma$
unbounded in $\GG(K_v)$, which is absurd.

We are now in a situation where $G\times H$ maps to $G\times Q$ with cocompact finite covolume image and
injectively on $\Gamma$; therefore the discreteness of the image of $\Gamma$ implies that this map is proper and
hence $H$ is a compact extension of $Q$. Pushing forward the measure on $(G\times H)/\Gamma$, we see that the
image of $\Gamma$ in $G\times Q$ is a lattice. Now $\Gamma$ is finitely generated (see above references
to~\cite[IX]{Margulis}) and thus the proof is completed as in Proposition~\ref{prop:SuperrigidPairs}. The
discrete factor occurring in the conclusion of the latter proposition is trivial for the same reason as in the
proof of Theorem~\ref{thm:arith:SxD}.
\end{proof}

\subsection{Lattices in general products}
We begin with the special case of totally disconnected groups.

\begin{proof}[Proof of Theorem~\ref{thm:SuperrigidPair:td}]
An issue that we need to deal with is that the projection of $\Gamma$ to $G_1$ is \emph{a priori }not
faithful. In order to circumvent this difficulty, we proceed to a preliminary construction.

Let $\iota: \Gamma \to \widehat \Gamma$ be the canonical map to the profinite completion of $\Gamma$ and denote its
kernel by $\Gamma\fire$; in other words, $\Gamma\fire$ is the finite residual of $\Gamma$\index{residual!finite}.
Let $\widehat G_1$ denote the locally compact group which is defined as the closure of the image of $\Gamma$ in
$G_1 \times \widehat \Gamma$ under the product map $\proj_1 \times \iota$, where $\proj_1 : G \to G_1$ is the
canonical projection. Since $\proj_1(\Gamma)$ is dense in $G_1$ and $\widehat \Gamma$ is compact, the canonical
map $\widehat G_1 \to G_1$ is surjective. In other words, the group $\widehat G_1$ is a compact extension of $G_1$.

We now define $G_1' = G_2 \times \cdots \times G_n$ and $\widehat G = \widehat G_1 \times G_1'$. Then $\Gamma$
admits a diagonal embedding into $\widehat G$ through which the injection of $\Gamma$ in $G$ factors. We will
henceforth identify $\Gamma$ with its image in $\widehat G$ and consider $\Gamma$ as an irreducible lattice of
$\widehat G$.

\smallskip \noindent
\emph{We claim that the pair $(\Gamma, \widehat G)$ is superrigid}.

\smallskip
The argument is a variation on the proof of Proposition~\ref{prop:commensurators}. Let $\tau : \Gamma \to
\HH(k)$ be as in the definition of superrigid pairs. Since $\tau(\Gamma)$ is finitely generated and linear, it
is residually finite~\cite{Malcev40}. This means that $\tau$ factors through $\ul\Gamma:=\Gamma/\Gamma\fire$.
Let $U < G'_1$ be a compact open subgroup, $\Gamma_U  = \Gamma \cap (\widehat G_1 \times U)$ and
$\ul{\Gamma_U}=\Gamma_U/(\Gamma_U\cap\Gamma\fire)$. By construction and Lemma~\ref{lem:CompactOpenTrick}, we can
consider $\ul{\Gamma_U}$ as a lattice in $\widehat G_1$ commensurated by $\ul\Gamma$. Arguing as in
Proposition~\ref{prop:commensurators}, when $\tau(\Gamma_U)$ is unbounded one applies commensurator
superrigidity yielding a continuous map $J \to \HH(k)$ and extending the map $\ul\Gamma\to\HH(k)$ and hence also
$\tau$. When $\tau(\Gamma_U)$ is bounded, one applies Lemma~\sref{lem:elementary:superrigidity} instead and the
resulting extension factors through $G'_1$. This proves the claim.

\smallskip
In order to conclude that the pair $(\Gamma, G)$ is also superrigid, it now suffices to apply Lemma~\ref{lem:containsH+}.
\end{proof}

\begin{cor}\label{cor:arith:td}
Theorem~\ref{thm:arith:general} holds in the particular case of totally disconnected groups.
\end{cor}

\begin{proof}
Theorem~\ref{thm:SuperrigidPair:td} provides the hypothesis needed for Proposition~\ref{prop:SuperrigidPairs}.
\end{proof}

We now turn to the general case $\Gamma<G=G_1\times \cdots \times G_n$ of Theorem~\ref{thm:arith:general}.
The main part of the remaining proof will consist of a careful analysis of how the lattice $\Gamma$ might sit
in various subproducts hidden in the factors $G_i$ or their finite index subgroups once the amenable radical has been trivialised.
It will turn out that $\Gamma$ is virtually a direct product $\Gamma'\times \Gamma''$, where $\Gamma'$ is an irreducible
lattice in a product $S'\times D'$ with $S'$ a semi-simple Lie (virtual) subproduct of $G$ and $D'$ a totally disconnected subgroup
of $G$ whose position will be clarified; as for $\Gamma''$, it is an irreducible lattice in a semi-simple Lie group $S''$ that turns
out to satisfy the assumptions of Margulis' arithmeticity. Of course, any of the above factors might well be trivial.

\begin{proof}[Completion of the proof of Theorem~\ref{thm:arith:general}]
The amenable radical is compact by Corollary~\ref{cor:radical:rigid} and hence we can assume that it is trivial.
The group $G$ (and hence also all finite index subgroups and factors) is compactly generated by
Lemma~\ref{lem:UGammaU}. Upon regrouping the last $n-1$ factors and in view of the definition of an irreducible
lattice (see p.~\pageref{def:IrreducibleLattice}), we can assume $G=G_1\times G_2$. We apply the solution to
Hilbert's fifth problem (compare Theorem~\sref{thm:folklore}) and write $G_i=S_i\times D_i$ after replacing $G$
and $\Gamma$ with finite index subgroups. Here $S_i$ are connected semi-simple centre-free Lie groups without
compact factors and $D_i$ totally disconnected compactly generated with trivial amenable radical. Set
$S=S_1\times S_2$ and $D=D_1\times D_2$. Thus $\Gamma$ is a lattice in $G = S \times D$. Notice that if $S$ is
trivial, then $G$ is totally disconnected and we are done by Theorem~\ref{thm:SuperrigidPair:td}. We assume
henceforth that $S$ is non-trivial. The main remaining obstacle is that the lattice $\Gamma$ need not be
irreducible with respect to the product decomposition $G = S \times D$.

\smallskip \noindent%
\emph{Observe that the closure $\overline{\proj_D(\Gamma)}$ of the projection of $\Gamma$ to $D$ has
trivial amenable radical. }

\smallskip%
Indeed $\proj_{D_i}(\Gamma)$ is dense in $\Gamma_i$ for $i= 1, 2$, hence the projection
$\overline{\proj_D(\Gamma)} \to D_i$ has dense image. The desired claim follows since $G$, and hence $D_i$, has
trivial amenable radical.

\smallskip
Let $U < D$ be a compact open subgroup and set $\Gamma_U = \Gamma \cap (S \times U)$. By
Lemma~\ref{lem:CompactOpenTrick}, the projection $\proj_S(\Gamma_U)$ of $\Gamma_U$ to $S$ is a lattice which is
commensurated by $\proj_S(\Gamma)$. The lattice $\proj_S(\Gamma_U)$ possesses a finite index subgroup which
admits a canonical splitting into finitely many irreducible groups $\Gamma^1 \times \cdots \times \Gamma^r$,
compare Theorem~\ref{thm:lattices:Irred}. Furthermore each $\Gamma^i$ is an irreducible lattice in a semi-simple
subgroup $S^i < S$ which is obtained by regrouping some of the simple factors of $S$.

Since the projection of $\Gamma$ to each $G_1$ and $G_2$, and hence to $S_1$ and $S_2$, is dense, it follows
that the projection of $\Gamma$ to each simple factor of $S$ is dense. We now consider the projection of
$\Gamma$ to the various factors $S^i$. In view of the preceding remark and the fact that $\Gamma^i$ is an
irreducible lattice in $S^i$, it follows that $\proj_{S^i}(\Gamma)$ is either dense in $S^i$ or discrete and
contains $\Gamma^i$ with finite index, see~\cite[IX.2.7]{Margulis}. Let now
$$S' = \la S^i \; | \; \proj_{S^i}(\Gamma) \text{ is non-discrete} \ra
\text{ and }%
S'' = \la S^i \; | \; \proj_{S^i}(\Gamma) \text{ is discrete} \ra.$$

\smallskip \noindent%
\emph{We claim that the projection of $\Gamma$ to $S'$ is dense. }

\smallskip%
If this failed, then by~\cite[IX.2.7]{Margulis} there would be a subproduct of some simple factors of $S'$ on
which the projection of $\Gamma $ is a lattice. Since each $\Gamma^i$ is irreducible, this subproduct is a
regrouping
$$S^{i_1} \times \cdots \times S^{i_p}$$
of some factors $S^i$. Now the projection of $\Gamma$ is a lattice in this subgroup, hence it contains the
product $\Gamma^{i_1} \times \cdots \times \Gamma^{i_p}$ with finite index and thus projects discretely to each
$S^{i_j}$. This contradicts the definition of $S'$ and proves the claim.

\smallskip \noindent%
\emph{Our next claim is that $\Gamma$ has a finite index subgroup which splits as $\Gamma' \times \Gamma''$,
where $\Gamma'' = \proj_{S''}(\Gamma)$ and $\Gamma'$ is a lattice in $S' \times D$.}

\smallskip
In order to establish this, we define
$$\Gamma' = \Ker(\proj: \Gamma \to S'') \hspace{1cm} \text{and} \hspace{1cm} \Gamma'' = \bigcap_{\gamma \in
\Gamma} \gamma\Gamma_U \gamma\inv.$$
Notice that $\Gamma'$ and $\Gamma''$ are both normal subgroups of $\Gamma$. Since
$\overline{\proj_{D}(\Gamma'')} $ is a compact subgroup of $D$ normalised by $\overline{\proj_D(\Gamma)}$, which
has trivial amenable radical, it follows that $\Gamma'' \subset S' \times S'' \times 1$. Therefore, the
intersection $\Gamma' \cap \Gamma''$ is a normal subgroup of $\Gamma$ contained in $S' \times 1 \times 1$. In
view of the preceding claim, we deduce that $\Gamma' \cap \Gamma'' = 1$. Thus $\la \Gamma' \cup \Gamma'' \ra
\lhd \Gamma$ is isomorphic to $\Gamma' \times \Gamma''$. Since $\Gamma'_U = \Gamma' \cap \Gamma_U$ projects to a
lattice in $S'$ which commutes with the projection of $\Gamma''$, we deduce moreover that
$\proj_{S'}(\Gamma'')=1$, or equivalently that $\Gamma'' < 1 \times S'' \times 1$.

Since the projection of $\Gamma$ to $S''$ has discrete image by definition, it follows from
Proposition~\ref{prop:Raghunathan} that $\Gamma' < S' \times 1 \times D$ projects onto a lattice in $S' \times
D$. On the other hand, the very definition of $S''$ implies $\proj_{S''}(\Gamma)$ contains
$\proj_{S''}(\Gamma_U)$, and hence also $\proj_{S''}(\Gamma'')$,  as a finite index subgroup. In particular,
this shows that $\Gamma' \times \Gamma''$ is a lattice in $S' \times S'' \times D$. Since it is contained in the
lattice $\Gamma$, we finally deduce that the index of $\Gamma' \times \Gamma''$ in $\Gamma$ is finite.

\smallskip
We observe that we have in particular obtained a lattice $\Gamma''<S''$ with $S''$ non-simple and $\Gamma''$ irreducible
(unless both $\Gamma''$ and $S''$ are trivial), because the projection of $\Gamma$ to any \emph{simple} Lie group factor
is dense: indeed, any simple factor must be a factor of some $G_i$ and $\Gamma$ projects densely on $G_i$.
It follows from Margulis' arithmeticity theorem~\cite[Theorem~1.(1')]{Margulis} that $\Gamma''$ is an arithmetic lattice
in $S''$.

\smallskip
Turning to the other lattice,
we remark that $\Gamma'$ admits a faithful Zariski-dense representation in a semi-simple group, obtained by
reducing the given representation of $\Gamma$. Furthermore, notice that the projection of $\Gamma$  to $S'$
coincides (virtually) with the projection of $\Gamma$. In particular it has dense image. Therefore, setting $D'
= \overline{\proj_D(\Gamma')}$, we may now view  $\Gamma' $ as an irreducible lattice in $S' \times D'$. We may
thus apply Theorem~\ref{thm:arith:SxD}. Notice that the same argument as before shows that $D'$ has trivial
amenable radical.

\smallskip \noindent%
\emph{We claim that the canonical discrete kernel $\Gamma'_{D'}$ is in fact a direct factor of $D'$. }

\smallskip
Indeed, since $\Gamma$ is residually finite by Malcev's theorem~\cite{Malcev40},
Proposition~\ref{prop:ResiduallyFinite} ensures that $\Gamma'_{D'}$ centralises the discrete residual
$D'^{(\infty)}$. In particular $ D'^{(\infty)} \cap \Gamma'_{D'}  = 1$ since $D'$ has trivial amenable radical.
Furthermore, since $D'/\Gamma'_{D'}$ is a semi-simple group, its discrete residual has finite index. In
particular, upon replacing $D'$ by a finite index subgroup we have $D' \cong  D'^{(\infty)} \times \Gamma'_{D'}$
as desired. It also follows that $\Gamma'_{D'}$ itself admits a Zariski-dense representation in a semi-simple
group.

\smallskip
It remains to consider again the projection maps $\proj_{D_i} : D \to D_i$. Restricting these maps to
$D'$ and using the fact that $\proj_{D_i}(D')$ is dense, we obtain that $D_i \cong \proj_{D_i}( D'^{(\infty)})
\times D'_i$, where $D'_i = \overline{\proj_{D_i}(\Gamma'_{D'})}$. The final conclusion follows by applying
Corollary~\ref{cor:arith:td} to the irreducible lattice $\Gamma'_{D'} < D'_1 \times D'_2$.
\end{proof}

\begin{proof}[Proof of Corollary~\ref{cor:arith:SimpleFactors}]
Since $\Gamma$ is finitely generated and irreducible, all $G_i$ are compactly generated (alternatively, apply
Lemma~\ref{lem:UGammaU}).

We claim that all projections $\Gamma\to G_i$ are injective. Indeed, if not, then (by induction on $n$) there is
$j$ such that the canonical discrete kernel $\Gamma_{G_j}$ is non-trivial. It is then cocompact, which implies
that the projection
$$ G_j/ \Gamma_{G_j} \times \prod_{i\neq j} G_i\ \lra\ \prod_{i\neq j} G_i$$
is proper. This is contradicts the fact that the image of $\Gamma$ in the left hand side above is discrete
whilst it is dense in the right hand side, proving the claim.

Suppose given a linear representation of $\Gamma$ in characteristic~$\neq 2,3$ whose image is not virtually
soluble. Arguing as in~\cite{Monod_alternative}, we can reduce to the case where we have a Zariski-dense
representation $\tau:\Gamma\to \HH(K)$ in a non-trivial connected adjoint absolutely simple group $\HH$ over a
finitely generated field $K$. Since $\tau(\Gamma)$ is infinite, we can choose a completion $k$ of $K$ for which
$\tau(\Gamma)$ is unbounded~\cite[2.1]{Breuillard-Gelander}.

Part of the argument in~\cite{Monod_alternative} is devoted to proving that the representation is \emph{a
posteriori} faithful. One can adapt the entire proof to the present setting, but we propose an alternative line
of reasoning using an amenability theorem from~\cite{Bader-Shalom}. Suppose towards a contradiction that the
kernel $\Gamma_0\lhd \Gamma$ of $\tau$ is non-trivial. Since the projections are injective, the closure $N_i$ of
the image of $\Gamma_0$ in $G_i$ is a non-trivial closed subgroup, which is normal by irreducibility and hence
is cocompact. Then Theorem~1.3 in~\cite{Bader-Shalom} implies that $\Gamma/\Gamma_0$ is amenable, contradicting
the fact that $\tau(\Gamma)$ is not virtually soluble in view of Tits' alternative~\cite{Tits72}.

At this point we can conclude by Theorem~\ref{thm:arith:general}.
\end{proof}

\section{Geometric arithmeticity}
\label{sec:geometric:arith}

\subsection{\cat lattices and parabolic isometries}
We now specialise the various arithmeticity results of Section~\ref{sec:arith} to the case of lattices in \cat
spaces and combine them with some of our geometric results.

Recall that a parabolic isometry is called \textbf{neutral}\index{isometry!neutral parabolic} if it has zero
translation length; the following contains Theorem~\ref{thm:Arithmeticity:neutral-intro} from the Introduction.

\begin{thm}\label{thm:ArithmeticityCAT0}
Let $X$ be a proper \cat space with cocompact isometry group and $\Gamma < G:=\Isom(X)$ be a finitely generated
lattice. Assume that $\Gamma$ is irreducible and that $G$ contains a neutral parabolic isometry. Then one of the
following assertions holds:
\begin{enumerate}
\item $G$ is a non-compact simple Lie group of rank one with trivial centre.

\item There is a subgroup $\Gamma_D \se \Gamma$ normalised by $G$, which is either finite or infinitely
generated and such that the quotient $\Gamma/\Gamma_D$ is an arithmetic lattice in a product of semi-simple Lie
and algebraic groups.
\end{enumerate}
\end{thm}

\begin{proof}
Let $X' \subseteq X$ be the canonical subspace provided by Theorem~\ref{thm:Lattice=>NoFixedPointAtInfty};
notice that $X'$ still admits a neutral parabolic isometry. Theorem~\sref{thm:Decomposition} and its addendum now
apply to $X'$. The space $X'$ has no Euclidean factor: indeed, otherwise
Theorem~\ref{thm:lattices:EuclideanSplitting} would imply $X'=\RR$, which has no parabolic isometries. The
kernel of the $\Gamma$-action on $X'$ is finite and we will include it in the subgroup $\Gamma_D$ below.

We distinguish two cases according as $X'$ has one or more factors.

In the first case, $\Isom(X')$ cannot be totally disconnected since otherwise Corollary~\sref{cor:BasicTotDisc}
point~\eqref{structure-pt:BasicTotDisc:elliptic}
rules out neutral parabolic isometries. Thus $\Isom(X')$ is a non-compact simple Lie group with trivial centre. If its real rank
is one, we are in case~(i); otherwise, $\Gamma$ is arithmetic by Margulis' arithmeticity theorem~\cite[Theorem~1.(1')]{Margulis}
and we are in case~(ii).

\smallskip

For the rest of the proof we treat the case of several factors for $X'$; let $\Gamma^*$ and let $H_{\Gamma^*}$ be as in
Section~\ref{sec:hull}. Note that $H_{\Gamma^*}$ acts cocompactly
on each irreducible factor of $X'$. Furthermore, each irreducible factor of $H_{\Gamma^*}$ is non-discrete
by Theorem~\ref{thm:lattices:Irred}.
Therefore $H_{\Gamma^*}$ is a product of the form $S \times D$ (possibly with one trivial factor), where $S$ is a
semi-simple Lie group with trivial centre and $D$ is a compactly generated totally disconnected group without discrete
factor.

By Corollary~\sref{cor:BasicTotDisc} point~\eqref{structure-pt:BasicTotDisc:elliptic}, the existence of a neutral parabolic
isometry in $G$ implies that $\Isom(X')$ is not totally disconnected.  Lemma~\ref{lem:hull} ensures that the
identity component of $\Isom(X')$ is in fact contained in $H_{\Gamma^*}$. Therefore, upon passing to a finite
index subgroup, the identity component of $\Isom(X')$ coincides with $S$.

If $D$ is trivial, then $H_{\Gamma^*} = S$ is a connected semi-simple Lie group containing $\Gamma$ as an
irreducible lattice. Since $S$ is non-simple, it has higher rank and we may appeal again to Margulis' arithmeticity theorem;
thus we are done in this case.

Otherwise, $D$ is non-trivial and we may then apply Theorem~\ref{thm:arith:SxD}. It remains to check that the
normal subgroup $\Gamma_D < \Gamma$, if non-trivial, is not finitely generated. But we know that $\Gamma_D$ is a
discrete normal subgroup of $D$. By Theorem~\ref{thm:density}, the lattice $\Gamma$, and hence also
$H_{\Gamma^*}$, acts minimally without fixed point at infinity on each irreducible factor of $X'$. Therefore,
Corollary~\sref{cor:NoEuclideanFactor} ensures that $D$ has no finitely generated discrete normal subgroup, as
desired.
\end{proof}

Here is another variation, of a more geometric flavour; this time, it is not required that there be a
\emph{neutral} parabolic isometry:

\begin{thm}\label{thm:arith:geometric}
Let $X$ be a proper geodesically complete \cat space with cocompact isometry group and $\Gamma < \Isom(X)$ be a
finitely generated lattice. Assume that $\Gamma$ is irreducible and residually finite.

If $G:=\Isom(X)$ contains any parabolic isometry, then $X$ is a product of symmetric spaces and Bruhat--Tits
buildings. In particular, $\Gamma$ is an arithmetic lattice unless $X$ is a real or complex hyperbolic space.
\end{thm}

\begin{proof}
We maintain the notation of the previous proof and follow the same arguments. We do not know a priori whether
there exists a neutral parabolic isometry. However, under the present assumption that $X$ is geodesically
complete, Corollary~\sref{cor:BasicTotDisc} point~\eqref{structure-pt:BasicTotDisc:completess} shows that the existence of
\emph{any} parabolic isometry is enough to ensure that $\Isom(X')$ is not totally disconnected.
Thus the conclusion of Theorem~\ref{thm:ArithmeticityCAT0} holds. In case~(i),
Theorem~\sref{thm:algebraic} point~\eqref{structure-pt:algebraic:homo} ensures that $X$ is a rank one symmetric space and we are
done. We now assume that (ii) holds and define $D$ as in the proof of Theorem~\ref{thm:ArithmeticityCAT0}.

The canonical discrete kernel $\Gamma_D$ is trivial by Theorem~\ref{thm:ResiduallyFiniteCAT0Lattice}. Since $D$
has no non-trivial compact normal subgroup by Corollary~\sref{cor:NoEuclideanFactor}, it follows from
Theorem~\ref{thm:arith:SxD} that $D$ is a totally disconnected semi-simple algebraic group. Therefore, the
desired result is a consequence of Theorem~\sref{thm:algebraic} point~\eqref{structure-pt:algebraic:homo}.
\end{proof}

For the record, we propose a variant of Theorem~\ref{thm:arith:geometric}:

\begin{thm}\label{thm:arith:geometric_bis}
Let $X$ be a proper geodesically complete \cat space with cocompact isometry group and $\Gamma < \Isom(X)$ be a
finitely generated lattice. Assume that $\Gamma$ is irreducible and that every normal subgroup of $\Gamma$ is
finitely generated.

If $G:=\Isom(X)$ contains any parabolic isometry, then $X$ is a product of symmetric spaces and Bruhat--Tits
buildings of total rank~$\geq 2$. In particular, $\Gamma$ is an arithmetic lattice.
\end{thm}

\begin{proof}
As for Theorem~\ref{thm:arith:geometric}, we can apply Theorem~\ref{thm:ArithmeticityCAT0}. We claim that case~(i)
is ruled out under the current assumptions. Indeed, a lattice in a simple Lie group of rank one is relatively
hyperbolic (see~\cite{Farb} or~\cite{Osin}) and as such has numerous infinitely generated normal subgroups
(and is even SQ-universal, see~\cite{Gromov87} or~\cite{Delzant96} for the hyperbolic case and~\cite{Arzhantseva-Minasyan-Osin}
for the general relative case).
In case~(ii) the discrete kernel $\Gamma_D$ is trivial and rank one is excluded as in case~(i) if the group is Archimedean;
if it is non-Archimedean, then there are no non-uniform finitely generated lattices (see~\cite{BassLubotzky})
and thus $\Gamma$ is again Gromov-hyperbolic which contradicts the assumption on normal, subgroups as before.
\end{proof}

We can now complete the proof of some results stated in the Introduction.

\begin{proof}[Proof of Theorem~\ref{thm:arith:para-intro}]
If $\Gamma$ is residually finite, then Theorem~\ref{thm:arith:geometric} yields the desired conclusion; it
therefore remains to consider the case where $\Gamma$ is not residually finite. We follow the beginning of the
proof of Theorem~\ref{thm:arith:geometric} until the invocation of
Theorem~\ref{thm:ResiduallyFiniteCAT0Lattice}, since the latter no longer applies. However, we still know that
there is a non-trivial Lie factor in $\Isom(X')$ and therefore we apply Theorem~\ref{thm:arith:SxD} in order to
obtain the desired conclusion about the lattice $\Gamma$. As for the symmetric space factor of the space, it is
provided by Theorem~\sref{thm:algebraic} point~\eqref{structure-pt:algebraic:homo}.
\end{proof}

\begin{proof}[Proof of Corollary~\ref{cor:para-intro}]
One implication is given by Theorem~\ref{thm:arith:para-intro}. For the converse, it suffices to recall that
unipotent elements exist in all semi-simple Lie groups of positive real rank.
\end{proof}

\subsection{Arithmeticity of linear \cat lattices}

We start by considering \cat lattices with a linear non-discrete linear commensurator:

\begin{thm}\label{thm:arith:commensurator:cat0}
Let $X$ be a proper geodesically complete  \cat space  with cocompact isometry group and $\Gamma < \Isom(X)$ be
a finitely generated lattice. Assume that $\Isom(X)$ possesses a finitely generated subgroup $\Lambda$
containing $\Gamma$ as a subgroup of infinite index, and commensurating $\Gamma$.

If $X$ is irreducible and $\Lambda$ possesses a faithful finite-dimensional linear representation (in
characteristic~$\neq 2, 3$), then $X$ is a symmetric space or a Bruhat--Tits building; in particular $\Gamma$ is
an arithmetic lattice.
\end{thm}

\begin{remark}\label{rem:Coxeter}
Several examples of irreducible \cat spaces $X$ of dimension $>1$ admitting a discrete cocompact group of
isometries with a non-discrete commensurator in $\Isom(X)$ have been constructed by F.~Haglund~\cite{Haglund}
and A.~Thomas~\cite{Thomas} (see also~\cite[Théorème~A]{Haglund08} and~\cite{BarnhillThomas}). In all cases that
space $X$ is endowed with walls; in particular $X$ is the union of two proper closed convex subspaces. This
implies in particular that $X$ is not a Euclidean building. Therefore,
Theorem~\ref{thm:arith:commensurator:cat0} has the following consequence: in the aforementioned examples of
Haglund and Thomas, \emph{either the commensurator of the lattice is nonlinear, or it is the union of a tower of
lattices}. In fact, as communicated to us by F.~Haglund, for most of these lattices the commensurator contains
elliptic elements of infinite order; this implies right away that the commensurator is not an ascending union of
lattices and, hence, it is nonlinear. Note on the other hand that it is already known that $\Isom(X)$ is mostly
nonlinear in these examples, since it contains closed subgroups isomorphic to the full automorphism group of
regular trees.
\end{remark}

\begin{proof}[Proof of Theorem~\ref{thm:arith:commensurator:cat0}]
Since $X$ is irreducible and the case $X = \RR$ satisfies the conclusions of the theorem, we assume henceforth
that $X$ has no Euclidean factor.

The $\Isom(X)$-action on $X$ is minimal by Lemma~\sref{lem:cocompact:minimal} and has no fixed point at infinity
by Corollary~\ref{cor:Lattice=>NoFixedPointAtInfty}. In particular, we can apply
Theorem~\sref{thm:GeodesicallyComplete}: either $\Isom(X)$ is totally disconnected or it is simple Lie group with
trivial centre and $X$ is the associated symmetric space. In the latter case, Margulis' arithmeticity theorem
finishes the proof. We assume henceforth that $\Isom(X)$ is totally disconnected.

Let $G$ denote the closure of $\Lambda$ in $\Isom(X)$. Note that $G$ acts minimally without fixed point at
infinity, since it contains a subgroup, namely $\Gamma$, which possesses these properties by
Theorem~\ref{thm:density}. In particular $G$ has trivial amenable radical by
Theorem~\sref{thm:geometric_simplicity} and thus the same holds for the dense subgroup $\Lambda < G$. In
particular any faithful representation of $\Lambda$ to an algebraic group yields a faithful representation of
$\Lambda$ to an adjoint semi-simple algebraic group with Zariski-dense image, to which we can apply
Theorem~\ref{thm:arith:commensurator}. As we saw, the group $G$ has no non-trivial compact (in fact amenable)
normal subgroup and furthermore $G$ is irreducible since $X$ is so, see Theorem~\sref{thm:geometric_simplicity}.
The fact that the lattice $\Gamma$ has infinite index in $\Lambda$ rules out the discrete case. Therefore $G$ is
a simple algebraic group and $\Gamma$ an arithmetic lattice.

It remains to deduce that $X$ has the desired geometric shape. This will follow from
Theorem~\sref{thm:algebraic}point~\eqref{structure-pt:algebraic:homo} provided we show that $\bd X$ is finite-dimensional and
that $G$ has full limit set. The first fact holds since $X$ is cocompact; the second is provided by
Corollary~\ref{cor:LimitSetLattice}.
\end{proof}

Remark~\ref{rem:Coxeter} illustrates that Theorem~\ref{thm:arith:commensurator:cat0} fails dramatically if one
assumes only that $\Gamma$ is linear. However, passing now to the case where $X$ is reducible, the linearity of
$\Gamma$ is enough to establish arithmeticity, independently of any assumption on commensurators,
the result announced in Theorem~\ref{thm:arith:geometric:lin-intro} in the Introduction.

\begin{thm}\label{thm:arith:geometric:lin}
Let $X$ be a proper geodesically complete \cat space with cocompact isometry group and $\Gamma < \Isom(X)$ be a
finitely generated lattice. Assume that $\Gamma$ is irreducible and possesses some faithful linear
representation (in characteristic~$\neq 2,3$).

If $X$ is reducible, then $\Gamma$ is an arithmetic lattice and $X$ is a product of symmetric spaces and
Bruhat--Tits buildings.
\end{thm}

\begin{proof}
In view of Theorem~\ref{thm:lattices:EuclideanSplitting}, we can assume that $X$ has no Euclidean factor. The
$\Isom(X)$-action on $X$ is minimal by Lemma~\sref{lem:cocompact:minimal} and has no fixed point at infinity by
Corollary~\ref{cor:Lattice=>NoFixedPointAtInfty}. In particular, we can apply
Theorem~\sref{thm:GeodesicallyComplete} to obtain decompositions of $\Isom(X)$ and $X$ in which the factors of
$X$ corresponding to connected factors of $\Isom(X)$ are isometric to symmetric spaces. There is no loss of
generality in assuming $\Gamma^* = \Gamma$ in the notation of Section~\ref{sec:hull}. Let now $G$ be the hull of
$\Gamma$. By Remark~\ref{rem:hull}, the group $\Gamma$ is an irreducible lattice in $G$.

Since $\Isom(X)$ acts minimally without fixed point at infinity, it follows from
Corollary~\ref{cor:LatticeNormaliser} that $\Gamma$ has trivial amenable radical. In particular any faithful
representation of $\Gamma$ to an algebraic group yields a faithful representation of $\Gamma$ to an adjoint
semi-simple algebraic group with Zariski-dense image, to which we can apply Theorem~\ref{thm:arith:general}.

The group $G$ has no non-trivial compact normal subgroup \emph{e.g.} by minimality. Furthermore the discrete
factor is trivial by Theorem~\ref{thm:lattices:Irred}. Therefore $G$ is a simple algebraic group and $\Gamma$ an
arithmetic lattice.

It remains to deduce that $X$ has the desired geometric shape and this follows exactly as in the proof of
Theorem~\ref{thm:arith:commensurator:cat0}.
\end{proof}

\subsection{A family of examples}\label{sec:example}
We shall now construct a family of lattices $\Gamma<G=S\times D$ as in the statement of
Theorem~\ref{thm:arith:SxD} (see also Theorem~\ref{thm:arith:geometric}) with the following additional
properties:

\begin{enumerate}
\item There is a proper \cat space $Y$ with $D<\Isom(Y)$ such that the $D$-action is cocompact, minimal and
without fixed point at infinity. In particular, setting $X=X_S\times Y$, where $X_S$ denotes the symmetric space
associated to $S$, the $\Gamma$-action on $X$ is properly discontinuous (in fact free), cocompact, minimal,
without fixed point at infinity.

\item The canonical discrete kernel $\Gamma_{\!D}\lhd D$ is infinite (in fact, it is a free group of countable
rank).

\item The profinite kernel of $D/\Gamma_{\!D} \to Q$ is non-trivial.
\end{enumerate}

\begin{remark}\label{rem:profinite_extension}
Since $D$ is minimal, it has no compact normal subgroup and thus we see that the profinite extension appearing
in Theorem~\ref{thm:arith:SxD} cannot be eliminated.
\end{remark}

\bigskip

We begin with a general construction:

\medskip

Let $\mathfrak g$ be (the geometric realisation of) a locally finite graph (not reduced to a single point) and
let $Q<\Isom(\mathfrak g)$ a closed subgroup whose action is vertex-transitive. In particular, $Q$ is a
compactly generated totally disconnected locally compact group. We point out that any compactly generated
totally disconnected locally compact group can be realised as acting on such a graph by considering Schreier
graphs $\mathfrak g$, see~\cite[\S\,11.3]{Monod_LN}; the kernel of this action is compact and arbitrary small. On
the other hand, if $Q$ is a non-Archimedean semi-simple group, one can also take very explicit graphs drawn on
the Bruhat--Tits building of $Q$, \emph{e.g.} the $1$-skeleton (this part is inspired by~\cite[1.8]{Burger-Mozes1},
see also~\cite{Burger-Mozes-Zimmer}).

Let moreover $C$ be an infinite profinite group and choose a locally finite rooted tree $\mathfrak t$ with a
level-transitive $C$-action for which every infinite ray has trivial stabiliser. For instance, one can choose
the coset tree associated to any nested sequence of open subgroups with trivial intersection, see the proof of
Théorème~15 in \S\,6 on p.~82 in~\cite{Serre_arbres}. We define a locally finite graph $\mathfrak h$ with a
$C\times Q$-action as the $1$-skeleton of the square complex $\mathfrak t\times \mathfrak g$. Let $\mathfrak
a=\wt{\mathfrak h}$ be the universal cover of $\mathfrak h$, $\Lambda=\pi_1(\mathfrak h)$ and define the totally
disconnected locally compact group $D$ by the corresponding extension
$$1\lra \Lambda\lra D \lra C\times Q\lra 1.$$

\begin{prop}\label{prop:construction}
There exists a proper \cat space $Y$ such that $D$ sits in $\Isom(Y)$ as a closed subgroup whose action is
cocompact, minimal and without fixed point at infinity.
\end{prop}

\begin{proof}
One verifies readily the following:

\begin{lem}\label{lem:renomalization}
Let $\mathfrak a$ be (the geometric realisation of) a locally finite simplicial tree and $D<\Isom(\mathfrak a)$
any subgroup. Let $x\in \mathfrak a$ be a vertex and let $Y$ be the completion of the metric space obtained by
assigning to each edge of $\mathfrak a$ the length $2^{-r}$, where $r$ is the combinatorial distance from this
edge to the nearest point of the orbit $D.x$.

Then $Y$ is a proper \cat space with a cocompact continuous isometric $D$-action. Moreover, if the $D$-action on
$\mathfrak a$ was minimal or without fixed point at infinity, then the corresponding statement holds for the
$D$-action on $Y$.\hfill\qedsymbol
\end{lem}

Apply the lemma to the tree $\mathfrak a=\wt{\mathfrak h}$ considered earlier. We claim that the $D$-action on
$\mathfrak a$ is minimal. Clearly it suffices to show that the $\Lambda$-action is minimal. Note that $\Lambda$
acts transitively on each fibre of $p:\wt{\mathfrak{h}}\to \mathfrak{h}$. Thus it is enough to show that the
convex hull of a given fibre meets every other fibre. Consider two distinct vertices $v, v' \in \mathfrak{h}$.
The product nature of $\mathfrak{h}$ makes it clear that $v$ and $v'$ are both contained in a common minimal
loop based at $v$. This loop lifts to a geodesic line in $\wt{ \mathfrak{h}}$ which meets the respective fibres
of $v$ and $v'$ alternatively and periodically. In particular, this construction yields a geodesic segment
joining two points in the fibre above $v$ and containing a point sitting above $v'$, whence the claim.

Since $\Lambda$ acts freely and minimally on the tree $\mathfrak{a}$ which is not reduced to a line, it follows
that $\Lambda$ fixes no end of $\mathfrak{a}$. Thus the lemma provides a proper \cat space $Y$ with a cocompact
minimal isometric $D$-action, without fixed point at infinity. It remains to show that $D < \Isom(Y)$ is closed.
This holds because the totally disconnected groups $\Isom(\mathfrak{a})$ and $ \Isom(Y)$ are isomorphic; indeed,
the canonical map $\mathfrak{a} \to Y$ induces a continuous surjective homomorphism $\Isom(\mathfrak{a}) \to
\Isom(Y)$, which is thus open.
\end{proof}

\begin{remark}\label{rem:non-smooth}
The above construction gives an example of a proper \cat space with a totally disconnected cocompact and minimal
group of isometries such that not all point stabilisers are open. Consider indeed the points added when
completing. Their stabilisers map to $Q$ under $D\to (C\times Q)$ and hence cannot be open. In other words, the
action is not \emph{smooth}\index{smooth} in the terminology of~\cite{CapraceTD}. Notice however that the set of
points with open stabiliser is necessarily a dense convex invariant set.
\end{remark}

We shall now specialise this general construction to yield our family of examples. Let $K$, $\HH$, $\Sigma$,
$K(\Sigma)$ be as described after Theorem~\ref{thm:arith:general} on~p.~\pageref{page:arithmetic_group}.
We write $\Sigma_f, \Sigma_\infty\se \Sigma$
for the subsets of finite/infinite places and assume that both are non-empty. Let $S=\prod_{v\in \Sigma_f}
\HH(K_v)^+$ and $Q=\prod_{v\in \Sigma_\infty} \HH(K_v)^+$. The group $\Delta=\HH(K(\Sigma))\cap(S\times Q)$ is
an irreducible cocompact lattice in $S\times Q$. Let $C$ be any profinite group with a dense inclusion
$\Delta\to C$. We now embed $\Delta$ diagonally in $S\times C\times Q$; clearly $\Delta$ is a cocompact lattice.
Let $\Gamma<G=S\times D$ be its pre-image. Then $\Gamma$ is a cocompact lattice since it contains the discrete
kernel of the canonical map $G \to S\times C\times Q$. It is clearly irreducible and therefore provides an
example that the structure of the description in the conclusion of Theorem~\ref{thm:arith:SxD} cannot be
simplified. Furthermore, the normal subgroup appearing in Theorem~\ref{thm:ArithmeticityCAT0}(ii) is also
unavoidable.

We end this section with a few supplementary remarks on the preceding construction:

\begin{enumerate}
\item If the profinite group $C$ has no discrete normal subgroup, then $\Gamma_{\! D} = \pi_1(\mathfrak{h})$
coincides with the quasi-centre of $D$. This would be the case for example if $C= \HH(K_v)$ and $\HH$ is almost
$K$-simple of higher rank, where $v$ is a non-Archimedean valuation such that $\HH$ is $K_v$-anisotropic. In
particular, in that situation $\Gamma_D$ is the unique maximal discrete normal subgroup of $D$ and the quotient
$D/\Gamma_D$ has a unique maximal compact normal subgroup. Thus the group $G$ admits a \emph{unique}
decomposition as in Theorem~\ref{thm:arith:SxD} in this case.

\item We emphasise that, even though $D/\Gamma_D$ decomposes as a direct product $C \times Q$ in the above
construction, the group $D$ admits no non-trivial direct product decomposition, since it acts minimally without
fixed point at infinity on a tree (see Theorem~\sref{thm:geometric_simplicity}).

\item The fact that $D/\Gamma_D$ decomposes as a direct product $C \times Q$ is not a coincidence. In fact, this
is happens always provided that every cocompact lattice in $S$ has the Congruence Subgroup Property (CSP).
Indeed, given a compact open subgroup $U$ of $D/\Gamma_D$, the intersection $\Gamma_U$ of $\Gamma \cap (S \times
U)$ is an irreducible lattice in $S \times U$ with trivial canonical discrete kernels. By (CSP), upon replacing
$\Gamma_U$ by a finite index subgroup (which amounts to replace $U$ by an open subgroup), the profinite
completion $\widehat{\Gamma_U}$ splits as the product over all primes $p$ of the pro-$p$ completions
$\widehat{(\Gamma_U)}_p$, which are just-infinite. Thus the canonical surjective map $\widehat{\Gamma_U} \to U$
shows that $U$ is a direct product. This implies that the maximal compact normal subgroup of $D/\Gamma_D$ is a
direct factor.

\item According to a conjecture of Serre's (footnote on page~489 in~\cite{Serre70}), if $S$ has higher rank then
every irreducible lattice in $S$ has (CSP). (See~\cite{Raghunathan04} for a recent survey on this conjecture.)
\end{enumerate}

\section{A few questions}

We conclude by collecting some further questions that we have encountered while working on this paper.

\medskip
It is well known that the Tits boundary of a proper \cat space with cocompact isometry group is
necessarily finite-dimensional (see~\cite[Theorem~C]{Kleiner}). It is quite possible that the same conclusion
holds under a much weaker assumption.

\begin{question}
Let $X$ be a proper \cat space such that $\Isom(X)$ has full limit set. Is the boundary $\bd X$
finite-dimensional?
\end{question}

A positive answer to this question would show in particular that the second set of assumptions~--- denoted~(b)~--- in
Theorem~\sref{thm:algebraic} is in fact redundant.

\medskip

Let $G$ be a simple Lie group acting continuously by isometries on a proper \cat space $X$.
The Karpelevich--Mostow theorem ensures that there exists a convex orbit when $X$ is a symmetric space of non-compact type.
This statement, however, cannot be generalised to arbitrary $X$ in view of Example~\sref{ex:parabolic_cone}.

\begin{question}\label{pb:KarpelevichMostow}
Let $G$ be a simple Lie group acting continuously by isometries on a proper \cat space $X$.
If the action is cocompact, does there exist a convex orbit?
\end{question}

It is shown in Theorem~\sref{thm:algebraic} point~\eqref{structure-pt:algebraic:homo} that if $X$ is geodesically
complete, then the answer is positive.
It is good to keep in mind Example~\sref{ex:Arbre-triangle},
which shows that the natural analogue of this question for a simple algebraic group over a non-Archimedean local
field has a negative answer.
More optimistically, one can ask for a convex orbit whenever the simple Lie group acts on a complete (not necessarily
proper) \cat space, but assuming the action non-evanescent (in the sense of~\cite{Monod_superrigid}). A positive answer
would imply superrigidity statements upon applying it to spaces of equivariant maps.

\medskip
Many of our statements on \cat lattices require the assumption of finite generation. One should of course wonder
for each of them whether it remains valid without this assumption. One instance where this question is
especially striking is the following (see Corollary~\ref{cor:Lattice=>NoFixedPointAtInfty}).

\begin{question}\label{pb:FiGen}
Let $X$ be a proper \cat space which is minimal and cocompact. Assume that $\Isom(X)$ contains a lattice. Is it
true that $\Isom(X)$ has no fixed point at infinity?
\end{question}

In a forthcoming article~\cite{Caprace-Monod_amen}, we shall establish a positive answer to this question by investigating
the rôle of unimodularity for the full isometry group.

\medskip
We have seen in Corollary~\sref{cor:NoOpenStabiliser} that if the isometry group of a proper \cat space $X$ is
non-discrete in a strong sense, then $\Isom(X)$ comes close to being a direct product of topologically \emph{simple}
groups.

\begin{question}
Retain the assumptions of  Corollary~\sref{cor:NoOpenStabiliser}. Is it true that $\soc(G^*)$ is a product of
simple groups? Is it cocompact in $G$, or at least does $G$ have compact Abelianisation?
\end{question}

Clearly Corollary~\sref{cor:NoOpenStabiliser} reduces the question to the case where $\Isom(X)$ is totally
disconnected. One can also ask if $\soc(G^*)$ is compactly generated (which is the case \emph{e.g} if it is cocompact in $G$).
If so, we obtain additional information by applying Proposition~\sref{prop:FilterNormalSubgroups}.

\medskip

In the above situation one furthermore expects that the geometry of $X$ is encoded in the structure of
$\Isom(X)$. In precise terms, we propose the following.

\begin{question}
Retain the assumptions of  Corollary~\sref{cor:NoOpenStabiliser}.
It it true that any proper cocompact action of $G$ on a proper \cat space $Y$ yields an
equivariant isometry $\bd X \to \bd Y$ between the Tits boundaries? Or an equivariant homeomorphism between the
boundaries with respect to the cône topology?
\end{question}

The discussion around Corollary~\ref{cor:arith:SimpleFactors} (\emph{cf.} Remark~\ref{rem:superrigidity:char})
suggests the following.

\begin{question}
Let $\Gamma<G=G_1\times\cdots\times G_n$ be an irreducible finitely generated lattice, where each $G_i$ is a
locally compact group. Does every character $\Gamma\to\RR$ extend continuously to $G$?
\end{question}

Y.~Shalom~\cite{Shalom00} proved that this is the case when $\Gamma$ is cocompact and in some other situations.

\def\cprime{$'$}
\providecommand{\bysame}{\leavevmode\hbox to3em{\hrulefill}\thinspace}


\end{document}